\RequirePackage{lineno}
\documentclass[journal]{IEEEtran}
\usepackage{times}
\usepackage{epsfig}
\usepackage{graphicx}
\usepackage{amsmath}
\usepackage{subfigure}
\usepackage[ruled,vlined]{algorithm2e}
\usepackage{amssymb}
\usepackage{multicol}
\usepackage{multirow}
\usepackage{enumitem}
\usepackage{stfloats}
\usepackage{amsthm}
\usepackage{cite}
\usepackage{color}
\usepackage{url}
\usepackage{relsize}
\usepackage{booktabs}
\usepackage{lipsum}
\usepackage{mathtools}
\usepackage{cuted}

\usepackage{url}
\usepackage{bm}
\usepackage{tikz} 
\usetikzlibrary{shapes.geometric, arrows}
\newcommand{\thickhline}{\noalign{\hrule height 1.0pt}}

\newtheorem{theorem}{Theorem}[]
\newtheorem{assumption}{Assumption}

\newcommand{\ten}[1]{\mathcal{#1}} 
\newcommand{\mat}[1]{\mathbf{#1}}
\newcommand{\vect}[1]{\boldsymbol{#1}}
\newcommand{\reff}[1]{(\ref{#1})}
\newcommand{\vecpar}{\boldsymbol{\xi}}

\newcommand{\multiGPC}{\Psi }
\newcommand{\multiGPCx}{\Phi }

\newcommand{\polyInd}{\alpha}
\newcommand{\basisInd}{\boldsymbol{\polyInd}}
\newcommand{\basisIndPar}{\boldsymbol{\beta}}

\newcommand{\yPC}{\sum\limits_{|\basisInd|=0}^{p} {c_{\basisInd}  \multiGPC_{\basisInd}  (\vecpar)} }

\newcommand{\zz}[1]{\textcolor{black}{#1}} 
\newcommand{\ccf}[1]{\textcolor{black}{#1}} 

\begin{document}

\title{Chance-Constrained and Yield-Aware  Optimization of Photonic ICs with Non-Gaussian Correlated Process Variations}

\author{Chunfeng Cui$^\star$, \and Kaikai Liu$^\star$, \and Zheng~Zhang,~\IEEEmembership{Member,~IEEE}
\thanks{$^\star$ C. Cui and  K. Liu contributed equally to this work. This work was partly supported by NSF Grant No. 1763699, NSF CAREER Award No. 1846476 and a UCSB start-up grant.}
\thanks{Chunfeng Cui, Kaikai Liu, and Zheng Zhang are in the Department of Electrical and Computer Engineering, University of California, Santa Barbara, CA 93106, USA (e-mail:  chunfengcui@ucsb.edu, kaikailiu@ucsb.edu, and zhengzhang@ece.ucsb.edu).}}


\maketitle
\begin{abstract}
    Uncertainty quantification has become an efficient tool for uncertainty-aware   prediction, but its power in yield-aware optimization has not been well explored from either theoretical or application perspectives. Yield optimization is a much more challenging task. On one side, optimizing the generally non-convex probability measure of  performance metrics is difficult. On the other side, evaluating the probability measure in each optimization iteration requires massive simulation data, especially when the process variations are non-Gaussian correlated. This paper proposes a data-efficient framework for the yield-aware optimization of photonic ICs. This framework optimizes the design performance with a yield guarantee, and it consists of two modules: a modeling module that builds stochastic surrogate models for design objectives and chance constraints with a few simulation samples, and a novel yield optimization module that handles probabilistic objectives and chance constraints in an efficient deterministic way. This deterministic treatment avoids repeatedly evaluating probability measures at each iteration, thus it only requires a few simulations in the whole optimization flow. We validate the accuracy and efficiency of the whole framework by a synthetic example and two photonic ICs. Our optimization method can achieve more than $30\times$ reduction of simulation cost and better design performance on the test cases compared with a Bayesian yield optimization approach developed recently.  
\end{abstract}


\IEEEpeerreviewmaketitle

\section{Introduction}
\IEEEPARstart{T}{he} demand for low-power, high-speed communications and computations have boosted the advances in  photonic integrated circuits. Based on the modern nano-fabrication technology, hundreds to thousands of photonic components can be integrated on a single chip \cite{nicholes2009world,kato200740}. However, process variations persist during all the fabrication processes and can cause a significant yield degradation in large-scale design and manufacturing \cite{chen2013process,lipka2016systematic,lu2017performance,pond2017predicting}. Photonic ICs are more sensitive to process variations (e.g., geometric uncertainties) due to their large device dimensions compared with the small wavelength. To achieve an acceptable yield, uncertainty-aware design optimization algorithms are highly desirable~\cite{weng2017stochastic}.  

Yield optimization algorithms try to increase the success ratio of a chip under random process variations, and they have been studied for a long time in the electronic circuit design 
\cite{yu1989efficient, li2011dram,wang2017efficient,wang2018efficient}.  However, it is still  expensive to reuse existing yield optimization solvers for photonic ICs. The major difficulties include: 1) the quantity of interest (e.g., the probability distribution of a bandwidth) does not admit an explicit expression. Instead, we only know the simulation values at  parameter sample points;  2) the design objectives and constraints are defined in a stochastic way. They are hard to compute directly and require massive numerical simulations to estimate their statistical distributions; 3) practical photonic IC designs often involve  non-Gaussian correlated process variations, which are more difficult to capture. To estimate the design yield efficiently, one alternative is to build a surrogate model. In~\cite{xu2005opera,li2006robust,li2007statistical}, posynomials were used to model statistical performance, and geometric programming was employed to optimize the worst-case performance.  The reference \cite{li2016efficient} proposed a  Chebyshev affine arithmetic method to predict the  cumulative distribution function. The recent Bayesian yield optimization \cite{wang2017efficient} approximated the probability density of the design variable under the condition of ``pass" by the kernel density estimation. The work 
\cite{wang2018efficient} further approximated the yield over the design variables directly by a Gaussian process regression. However, these machine learning techniques may still require many simulation samples. Furthermore, worst-case optimization or only optimizing the yield can lead to non-optimal (and even poor) chip performance.  

Recently, uncertainty quantification methods based on generalized polynomial chaos have
achieved great success in modeling the impact caused by process variations in electronic and photonic ICs  \cite{zzhang:tcad2013,zhang2014calculation,zzhang:huq_tcad,zhang2017big,manfredi2013uncertainty, vrudhula2006hermite, strunz2008stochastic,rufuie2013generalized,tao2007stochastic,shen2009variational, waqas2018uncertainty,waqas2018stochastic}. A novel stochastic collocation approach was further proposed in \cite{cui2018stochastic,  cui2018stochastic2} to handle non-Gaussian correlated process variations, which shows significantly better accuracy and efficiency than~\cite{weng2015uncertainty} due to an optimization-based quadrature rule. These techniques construct stochastic surrogate models with a small number of simulation samples, but their power in yield optimization has not been well explored despite  the recent robust optimization methods~\cite{weng2017stochastic} based on generalized polynomial chaos.

  Leveraging the chance-constrained optimization~\cite{li2008chance} and our recently proposed uncertainty quantification solvers~\cite{cui2018stochastic,  cui2018stochastic2}, this paper presents a data-efficient technique to optimize photonic ICs with non-Gaussian correlated process variations. Instead of just optimizing the yield, we optimize a target performance metric while enforcing the probability of violating  design rules to be smaller than a user-defined threshold. Doing so can avoid performance degradation in yield optimization. Chance-constrained optimization~\cite{li2008chance} has been widely used in system control~\cite{mesbah2014stochastic}, autonomous vehicles \cite{blackmore2010probabilistic}, and reliable power generation \cite{akhavan2018energy,wang2017chance}, but it has not been investigated for yield optimization of electronic or photonic ICs. Our specific contributions include:
\begin{itemize}
    \item  A chance-constraint optimization framework that can achieve high chip performance and high yield simultaneously under non-Gaussian correlated process variations. 
    
     \item A surrogate model that approximates the stochastic objective and constraint functions with a few simulations. Since both the objective function and constraints are only available through a  black-box simulator, we build a surrogate model based on the recent uncertainty quantification solver~\cite{cui2018stochastic2}. 
     The main step is to compute a quadrature rule in the joint space of  design variables and stochastic parameters by a new three-stage optimization process. 
    \item A deterministic reformulation. A major challenge of chance-constrained optimization is to reformulate the stochastic constraints into deterministic ones \cite{nemirovski2006convex}.  We reformulate the probabilistic objective function and constraints as 
    non-smooth deterministic functions. Afterward, we transform  \ccf{them} into an equivalent polynomial optimization, which can be  solved efficiently.
    \item Validations on benchmarks. Finally, we validate the efficiency of our proposed framework on a synthetic example, a microring add-drop filter, and a Mach-Zehnder filter. Preliminary numerical experiments show that our proposed framework can find the optimal design variable efficiently. 
    Compared with the Bayesian yield optimization method \cite{wang2017efficient}, our proposed method can reduce the number of simulations by $30\times$, achieve better performance,  and produce a similar yield  on the test cases.
\end{itemize}
This work should be regarded as a preliminary result in this direction, and many topics can be investigated in the future.





\section{Preliminaries}

\subsection{The Yield Optimization} 
The   yield  is defined as the percentage of qualified  products overall. 
For a photonic IC, denote the design variables by $\vect   x=[x_1,x_2,...,x_{d_1}]^T\in\mathcal{X}$ and the process variations by random parameters  $\vecpar=[\xi_1,\xi_2,...\xi_{d_2}]^T\in\Omega$. Suppose $\vect x$ is uniformly distributed in a bound domain and $\vecpar$  follows a  probability distribution $\rho(\vecpar)$. 
Let $\{y_i(\vect x, \vecpar)\}_{i=1}^n$ denote a set of performance metrics of interest,
$\{u_i\}_{i=1}^n$ denote its required upper bound, 
and  $I(\vect x, \vecpar)$  denote the indicator function:
\begin{equation}
    I(\vect x, \vecpar)=\left\{
    \begin{array}{cl}
        1, & \text{if   $y_i(\vect x, \vecpar)\ccf{\le u_i, \forall i=1,\ldots,n;}$} \\
        0, & \text{otherwise.}
    \end{array}
    \right.
    \label{equ:pass}
\end{equation}
\ccf{The yield at a certain design choice $\vect x$ is defined as 
\begin{equation}
    Y(\vect x)={\rm Prob}_{\vecpar}(\mat y(\vect x,\vecpar)\le\vect u | \vect x)=\mathbb{E}_{\vecpar}[I(\boldsymbol{x},\boldsymbol{\xi})].
\end{equation}
}
The yield optimization problem aims to find an optimal design variable $\vect{x}^*$ such that 
\begin{equation}\label{equ:max_yield}
    \boldsymbol{x}^*= \underset{\vect x\in \mathcal{X}}{\text{argmax}} \ \ \ccf{Y(\vect x)}. 
\end{equation}
 There are three major difficulties in solving the above yield optimization problem:   1) the indicator function $I(\vect x,\vecpar)$ does not always admit an explicit formulation; 2) computing the yield $\ccf{Y(\vect x)}$ involves  a  non-trivial numerical integration, which  requires  numerous simulations at each design variable $\vect x$; 3)   $\ccf{Y(\vect x)}$ is an implicit non-convex function and it is difficult to compute its optimal solution.

\subsection{Chance Constraints}
The chance constraint is a powerful technique in uncertainty-aware optimization \cite{li2008chance}. In comparison with the deterministic constraints or the worst-case constraints where the risk level $\epsilon$ is zero, a chance constraint enforces the probability of satisfying a stochastic constraint to be above a certain confidence level $1-\epsilon$ ($\epsilon$ is usually not zero): 
\begin{equation}\label{equ:chancecons}
    \text{Prob}_{\vecpar}(y(\boldsymbol{x},\boldsymbol{\xi})\le \vect u)\ge1-\epsilon 
\end{equation}
or equivalently, the probability of violating the constraint to be smaller than the risk level $\epsilon$:
\begin{equation}\label{equ:chancecons2}
    \text{Prob}_{\vecpar}(y(\boldsymbol{x},\boldsymbol{\xi})\ge \vect u)\le\epsilon.
\end{equation}
 
Under  strict conditions, such as the parameters being independent and  $y(\vect x, \vecpar)$ being a linear function,   \reff{equ:chancecons} can be  reformulated into  equivalent deterministic constraints \cite{calafiore2005distributionally}. In other words, one can reformulate the left-hand side of \reff{equ:chancecons} by its  probability density function    and substitute the right-hand side by a constant related to the cumulative density function. However, these conditions rarely hold in practice. Even if the conditions hold, computing the probability density function or cumulative density function of an uncertain variable can be intractable \cite{zhang2014calculation,nemirovski2006convex}. In these cases, we seek for  deterministic reformulations that can well approximate the chance constraints. There is a trade-off in choosing the reformulation:  if the reformulation is aggressive (the feasible domain is enlarged),  it may result in an infeasible solution; Otherwise, if the reformulation is conservative (the feasible domain is decreased),  the solution may be degraded. 

One may convert the chance constraint  \reff{equ:chancecons} to a deterministic constraint via the mean and variance of $y(\boldsymbol{x},\boldsymbol{\xi})$ \cite{nemirovski2006convex,calafiore2005distributionally}: 
\begin{equation}\label{equ:chance2pop}
   \mathbb{E}_{\vecpar}[y(\boldsymbol{x},\boldsymbol{\xi})] + \kappa_\epsilon\sqrt{\text{var}_{\vecpar}[y(\boldsymbol{x},\boldsymbol{\xi})]} \ge \vect u.
\end{equation}
Here $\mathbb{E}_{\vecpar}[\cdot]$ denotes the mean value, $\text{var}_{\vecpar}[\cdot]$  denotes the variance. The constant $\kappa_\epsilon$ is chosen as 
$\kappa_\epsilon=\sqrt{(1-\epsilon)/\epsilon}$. The detailed proof is shown in Appendix~\ref{app:chance_proof}. 
 It is worth noting that \reff{equ:chance2pop} is a stronger constraint than \eqref{equ:chancecons}:  every feasible point of \reff{equ:chance2pop} is also a feasible point of the original chance constraint \reff{equ:chancecons}.

\subsection{Stochastic Spectral Methods}
Assume that $y(\vecpar)$  is a smooth function satisfying  $\mathbb{E}[y^2(\vecpar)]\le\infty$. The stochastic spectral methods can approximate $y(\vecpar)$ by   orthonormal polynomial basis functions:
\begin{equation}
\label{eq:ygpc}
y(\vecpar) \approx \yPC, \; {\rm with}\; \mathbb{E}\left[{\multiGPC}_{\basisInd} (\vecpar)\multiGPC_{\boldsymbol{\beta }}\left( \vecpar \right)\right ]=\delta_{\basisInd, \boldsymbol{\beta }}.
\end{equation}
Here $|\basisInd|=\alpha_1+\ldots+\alpha_{d_2}$, $\multiGPC_{\basisInd} (\vecpar)$ is an orthonormal basis function indexed by $\basisInd$, and $c_{\basisInd}$ is its corresponding coefficient. 

If  the parameters $\vecpar$ are independent,   $\rho(\vecpar)$ equals  the products of   its one-dimensional marginal density function $\rho_i(\xi_i)$. In this case, the basis function $\multiGPC_{\basisInd} (\vecpar)$ is the product of multiple one-dimensional orthogonal basis functions 
\begin{equation}
   \multiGPC_{\basisInd} (\vecpar) = \psi_1(\xi_1)\ldots\psi_{d_2}(\xi_{d_2}). 
\end{equation}
These one-dimensional basis functions $\psi_i(\xi_i)$ can be  constructed by the three term recursion \cite{Walter:1982}. Various stochastic spectral approaches have been proposed to compute the coefficients $c_{\basisInd}$, including the intrusive (i.e., non-sampling) solvers (e.g., stochastic Galerkin~\cite{sfem}, the stochastic testing~\cite{zzhang:tcad2013}) and the non-intrusive (i.e., sampling) solvers (e.g., stochastic collocation~\cite{col:2005}).  In the past few years, there has also been a rapid progress in handling  high-dimensional parameters, such as the tensor recovery method~\cite{zhang2017big}, the compressive sensing technique\cite{li2010finding}, ANOVA (analysis of variance) or HDMR (the high-dimensional model representation)~\cite{zzhang_cicc2014}, and the hierarchical uncertainty quantification~\cite{zzhang:huq_tcad}.

In practice, the random parameters may be correlated. 
 If the parameters $\vecpar$ are  non-Gaussian correlated, the computation is more difficult.  In such cases,   ${\multiGPC}_{\basisInd}(\vecpar)$ can be constructed by the Gram-Schmidt approach  \cite{cui2018stochastic, cui2018stochastic2} or the Cholesky factorization   \cite{cui2019high, cui2018uncertainty}. The main difficulty lies in computing  high order moments of $\vecpar$, which can be well resolved by the functional tensor train approach   \cite{cui2018uncertainty}.


\section{Our Yield-aware Optimization Model}

In this section, we  show our yield-aware chance constrained  optimization model, and  illustrate how to convert the stochastic formulation to a deterministic one.The basic assumptions are listed as follows. 

\begin{assumption}

We made the following assumptions:
\begin{enumerate}
    \item The design variable $\vect x$ is upper and lower bounded, i.e.,  $\vect x \in\ten X=[\mat a,\mat b]^{d_1}$; 
    
    \item The stochastic parameter $\vecpar \in\Omega\in\mathbb{R}^{d_2}$ admits a non-Gaussian correlated density function $\rho(\vecpar)$; 
    
    \item The yield is qualified by the following constraints:
   \begin{equation}\label{equ:success}
    y_i(\vect x,\vecpar)\le  u_i,\  \forall\, i\in[n].
\end{equation}
    Here $[n]=1,\ldots,n$ and $\mathbb{E}[y_i(\vect x,\vecpar)]\le u_i$. Each individual quantity   $y_i(\vect x,\vecpar)$ is a black-box function, and we can obtain its function values at given sample points.  
\end{enumerate}
\end{assumption}

\ccf{The design variables $\vect{x}$  are deterministic  without any probability measures, and all samples of $\vect{x}$ are equally important in the optimization process. Therefore, we treat $\vect{x}$ as  mutually independent random variables with a uniform distribution and use   Legendre polynomials as their basis functions. The process variations $\vecpar$ are non-Gaussian correlated, which enables our model to handle generic cases.
}
\subsection{The Probabilistic Yield Optimization Model}
The yield at a given design variable $\vect x$ can be defined as  the probability that the yield conditions \reff{equ:success} are satisfied, i.e.,  
\begin{equation*}
   Y(\boldsymbol{x})=\text{Prob}_{\vecpar}(  \vect y(\vect x,\vecpar)\le \vect u). 
\end{equation*}
Here, $\vect y(\vect x,\vecpar)=[y_1(\vect x,\vecpar),\ldots, y_n(\vect x,\vecpar)]^T$ and $\vect u=[u_1,\ldots,u_n]^T$. Consequently, the yield optimization problem can be described as: 
\begin{align}\label{equ:P00}
    \max_{\vect x\in\ten X}\quad&\text{Prob}_{\vecpar}(  \vect y(\vect x,\vecpar)\le \vect u).
\end{align}

However, the above yield maximization often contradicts with our performance goals. For instance, one may have to reduce the clock rate of a processor significantly in order to achieve a high yield. As a result, directly optimizing the yield may lead to an over-conservative design. 
\ccf{In practice, the design variables that provide the best yield may be nonunique, and we hope to chose a design that achieves good performance and high yield simultaneously.
Therefore, we ensure the yield with a chance constraint} 
\begin{equation}\label{equ:joint_chance}
     \text{Prob}_{\vecpar}(  \vect y(\vect x,\vecpar)\le \vect u)\ge 1-\epsilon,
\end{equation}
and   optimize the expected value of an uncertain   performance metric $f(\vect x,\vecpar)$ by the \ccf{following yield-aware optimization}:
\begin{align}\label{equ:P0_0}
\nonumber \min_{\vect x\in\ten X}\quad&\mathbb{E}_{\vecpar}[f(\vect x,\vecpar)]\\
\text{s.t.}\quad & \text{Prob}_{\vecpar}(  \vect y(\vect x,\vecpar)\le \vect u)\ge 1-\epsilon.
\end{align}
Here $\epsilon$ is a given risk level to control the yield. 
\ccf{The above formulation is not equivalent to (\ref{equ:P00}).} It can  describe, for instance, the following design optimization problem: minimize the average power consumption of a photonic IC while ensuring at least $95\%$ yield (i.e., with $5\%$ probability of violating timing and bandwidth constraints) under process variations. Note that $f(\vect x,\vecpar)$ may also be the function (e.g., weighted sum) of several performance metrics that we intend to optimize simultaneously.
\ccf{The parameter $\epsilon$ can help designers balance between the yield and a target performance goal (i.e., power consumption). A small $\epsilon$ results in a higher yield but possibly a worse performance metric. Therefore, the value of $\epsilon$ can be chosen adaptively and case-dependently by the users based on on their specific requirements on the performance and yield.}

Because the yield function $Y(\vect x)$ and the objective function $f(\vect x,\vecpar)$ are not available, we have to estimate the yield and objective at a certain design variable  $\vect x$  by the Monte Carlo method \cite{yu1989efficient,  li2011dram}. This requires a huge number of simulation samples at each design variable $\vect x$, which is infeasible for many simulation-expensive photonic IC design problems.

\begin{figure}
    \centering
    \begin{tikzpicture}
    
    \tikzstyle{process} = [rectangle, minimum width=3cm, minimum height=0.1cm, text centered, draw=black,text width=3cm]

 \node (pro1) [process] {\footnotesize Build the  chance constrained   model \reff{equ:P0_0}};
\node (pro2) [process, below of= pro1, yshift=-0.2cm] {\footnotesize Reformulate \reff{equ:P0_0} into   \reff{equ:P0} with $n$  constraints};
\node (pro3) [process, below of= pro2, yshift=-0.2cm] {\footnotesize Reformulate \reff{equ:P0} into   deterministic model  \reff{equ:deterministic_model}};
\node (pro4) [process, below of= pro3, yshift=-0.2cm] {\footnotesize Derive the polynomial optimization model \reff{equ:polyopt}};
\node (pro5) [process, below of= pro4, yshift=-0.2cm] {\footnotesize \ccf{Solve}  \reff{equ:polyopt} and output the optimal design}; 

\node (pro22) [process, right of= pro1, xshift=3cm] {\footnotesize Input the range of $\vect x$ and the PDF of $\vecpar$};
\node (pro23) [process, right of= pro2, xshift=3cm] {\footnotesize Solve  \reff{equ:co-optimQuad} to compute the quadrature rule};
\node (pro24) [process, right of= pro3, xshift=3cm] { \footnotesize Call the simulator at the   quadrature points};
\node (pro25) [process, right of= pro4, xshift=3cm] {\footnotesize Construct the surrogate model by \reff{equ:coeffs}};

\draw[->] (pro1.south) -- (pro2.north);
\draw[->] (pro2.south) -- (pro3.north);
\draw[->] (pro3.south) -- (pro4.north);
\draw[->] (pro4.south) -- (pro5.north);
 
\draw[->] (pro22.south) -- (pro23.north);
\draw[->] (pro23.south) -- (pro24.north);
\draw[->] (pro24.south) -- (pro25.north);

\draw[->] (pro25.west) -- (pro4.east); 
    
    \end{tikzpicture}
    \caption{The flowchart of our proposed framework for solving the chance constrained yield-aware  optimization.}
    \label{fig:flowchart}
\end{figure}
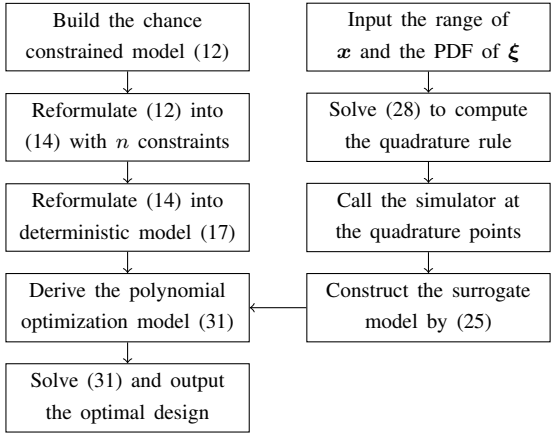



\ccf{
 Due to the ease of implementation, we reformulate the joint chance constraint in (\ref{equ:joint_chance}) into individual chance constraints: 
\begin{equation}\label{equ:yield2chance}
    \text{Prob}_{\vecpar}(y_i(\vect x,\vecpar)\le  u_i)\ge 1-\epsilon_i, \forall\, i\in[n].
\end{equation}
In this formulation, $\epsilon_i$ means the  risk tolerance of violating the $i$-th design specification. 
Since $\text{Prob}_{\vecpar}(\mat y(\vect x,\vecpar)\le\mat u)=\text{Prob}_{\vecpar}(\cap_{i=1}^n  (y_i(\vect x,\vecpar)\le u_i))=1-\text{Prob}_{\vecpar}(\cup_{i=1}^n   (y_i(\vect x,\vecpar)\ge u_i))$,  the probability of the joint chance constraint can be upper and lower bounded by the  individual chance constraints:
\begin{align*}
   \max_{i=1,\ldots,n}\ \text{Prob}_{\vecpar}(y_i(\vect x,\vecpar)\ge u_i) &
   \le  \text{Prob}_{\vecpar}(\cup_{i=1}^n (y_i(\vect x,\vecpar)\ge u_i)) \\
   & \le  \sum_{i=1}^n \text{Prob}_{\vecpar}(  y_i(\vect x,\vecpar)\ge u_i). 
\end{align*}
When $\epsilon_i=\epsilon$ for all $i$, (\ref{equ:yield2chance}) is a relaxation of  (\ref{equ:joint_chance}) (e.g., the feasible domain is enlarged); 
when $\sum_{i=1}^n\epsilon_i\le \epsilon$, (\ref{equ:yield2chance}) becomes more  conservative than  (\ref{equ:joint_chance}) (e.g., the feasible domain is reduced).   
In this paper, we do not give the universal best choice of $\epsilon_i$. 
Instead, the users can tune the parameters adaptively based on their requirements. 
}  



Consequently, we have the following chance-constrained yield-aware optimization model
\begin{align}\label{equ:P0}
\nonumber \min_{\vect x\in\ten X}\quad&\mathbb{E}_{\vecpar}[f(\vect x,\vecpar)]\\
\text{s.t.}\quad &  \text{Prob}_{\vecpar}(y_i(\vect x,\vecpar)\le  u_i)\ge 1-\epsilon_i, \forall\, i\in[n].
\end{align}





\subsection{The Deterministic  Reformulation}

The chance-constrained optimization problem \reff{equ:P0} is difficult to solve directly. This problem is more challenging when $y_i(\vect x, \vecpar)$ is nonlinear because it is almost impossible to formulate the chance constraints in \reff{equ:P0} to equivalent  deterministic formulations. 
A \ccf{naive} approach is to replace the stochastic constraints by  inequality constraints over the expected constraints:
\begin{align}\label{equ:P0_mean}
\nonumber \min_{\vect x\in\ten X}\quad&\mathbb{E}_{\vecpar}[f(\vect x,\vecpar)]\\
\text{s.t.}\quad &
    \mathbb{E}_{\vecpar}[y_i(\vect x,\vecpar)]\le  u_i, \forall\, i\in[n].
\end{align}
However, this treatment will lose the probability density  information and may not provide a high-quality solution, although it can help improve the yield in practice. We will illustrate this phenomenon in numerical experiments in Section~\ref{sec:synthetic}. 

\zz{Therefore, we do not use the formulation in \eqref{equ:P0_mean}.} Instead, we adopt the second-order moment approach   in  \cite{nemirovski2006convex,calafiore2005distributionally} and replace \reff{equ:yield2chance} by  
\begin{equation}\label{equ:determin_chance}
    \mathbb{E}_{\vecpar}[y_i(\vect{x},\vecpar)]+\kappa_{\epsilon_i} \sqrt{\text{var}_{\vecpar}[y_i(\vect{x},\vecpar)]}\le u_i, \forall\, i\in[n]. 
\end{equation}
Here,  $\kappa_{\epsilon_i}=\sqrt\frac{1-\epsilon_i}{\epsilon_i}$ is a scaling parameter. We present the detailed proof  in Appendix \ref{app:chance_proof} and 
  point out the following: 
\begin{itemize}
    \item Constraint \reff{equ:determin_chance} is a stronger condition than \reff{equ:yield2chance}. In other words, each feasible point of  \reff{equ:determin_chance} is also  feasible for the chance constraint \reff{equ:yield2chance}; 
    
    \item The parameter $\epsilon_i$ is a user-defined risk tolerance. When $\epsilon_i$ decreases, the feasible set will become smaller. However,  the optimal solution may result in a higher yield;  
    \item When the variance $\text{var}_{\vecpar}[y_i(\vect{x},\vecpar)]$ is small enough, the feasible set of \reff{equ:determin_chance} is  close to the deterministic constraint  $\mathbb{E}_{\vecpar}[y_i(\vect{x},\vecpar)]\le u_i$. 
\end{itemize}

Consequently, the probabilistic optimization model \reff{equ:P0} is reformulated into a deterministic optimization problem: 
  \begin{align}\label{equ:deterministic_model}
   \nonumber \min_{\vect x\in\ten X}\quad& \mathbb{E}_{\vecpar}[ f(\vect x,\vecpar)]\\
    \text{s.t.}\quad &\mathbb{E}_{\vecpar}[y_i(\vect{x},\vecpar)]+\kappa_{\epsilon_i} \sqrt{\text{var}_{\vecpar}[y_i(\vect{x},\vecpar)]}\le u_i, \forall\, i\in[n]. 
\end{align}



\section{Algorithm and Implementation Details}

We cannot solve problem \reff{equ:deterministic_model} directly because we do not know   the mean values and variances for the black-box functions $\{y_i(\vect x,\vecpar)\}_{i=1}^n$ and $f(\vect x,\vecpar)$.  A direct approach is to apply a Monte Carlo method to estimate the mean values and variances for every iterate $\vect x$. However, this is not affordable because of the large number of numerical simulations.

In this section,  we build the surrogate model for $f(\vect x, \vecpar)$ and  $\{y_i(\vect x,\vecpar)\}_{i=1}^n$ by using generalized polynomial chaos \cite{xiu2002wiener} and our recent developed uncertainty quantification solver \cite{cui2018stochastic, cui2018stochastic2}. Once the surrogate models are constructed, we can perform deterministic optimization. The main task is to build the orthogonal basis functions $\multiGPCx_{\basisInd}(\vect x)$ and $\multiGPC_{\basisIndPar}(\vecpar)$, and compute the coefficients $c_{\basisInd,\basisIndPar}^i$ and $h_{\basisInd,\basisIndPar}$ such that 
\begin{equation}\label{equ:surrogate}
    y_i(\vect x,\vecpar)\approx\sum_{|\basisInd|+|\basisIndPar|=0}^{p} c^i_{\basisInd,\basisIndPar}\multiGPCx_{\basisInd}(\vect x)\multiGPC_{\basisIndPar}(\vecpar),
\end{equation}
and \begin{equation}\label{equ:surrogate_obj}
    f(\vect x,\vecpar)\approx\sum_{|\basisInd|+|\basisIndPar|=0}^{p} h_{\basisInd,\basisIndPar}\multiGPCx_{\basisInd}(\vect x)\multiGPC_{\basisIndPar}(\vecpar).
\end{equation} 
Once the above surrogate models are obtained, the mean value of $y_i(\vect x,\vecpar)$ can be approximated by  
\begin{equation}\label{equ:mean}
    \mathbb{E}_{\vecpar}[y_i(\vect x,\vecpar)] \approx  \sum_{|\basisInd|=0}^{p}c^i_{\basisInd,\mat 0}\multiGPCx_{\basisInd}(\vect x),
\end{equation}
and the variance is approximated by
\begin{equation}\label{equ:var}
    \text{var}_{\vecpar}[y_i(\vect x,\vecpar)]  \approx \sum_{|\basisIndPar|=1}^{p}\left(\sum_{|\basisInd|=0}^{p-|\basisIndPar|} c^i_{\basisInd,\basisIndPar}\multiGPCx_{\basisInd}(\vect x)\right)^2.
\end{equation}
\ccf{Equation (\ref{equ:var}) is obtained based on the orthonormal property of the basis functions. The detailed proof is shown in Appendix B.} 
The mean value  of the objective function $f(\vect x,\vecpar)$ can be evaluated in the same way. Finally, the deterministic yield optimization model \reff{equ:deterministic_model} has an explicit expression and can be solved. 

The overall framework is summarized in Algorithm~\ref{alg:framework}. In the following, we explain the implementation details.

\begin{algorithm}[t]
\label{alg:framework}
\caption{Our Proposed Chance-Constrained Yield-aware Optimization Solver}
      \SetKwInput{Input}{Input}
      \SetKwInput{Output}{Output}
\Input{The range of the design variable $\vect x$,  probability density function of the non-Gaussian correlated random parameters $\rho(\vecpar)$, the polynomial order  $p$, the upper bounds of performance metrics $\{u_i\}_{i=1}^n$, and the chance constraint thresholds $\{\epsilon_i\}_{i=1}^n$.}
\begin{enumerate}[leftmargin=*]
    \item[1.] Construct the basis functions $\multiGPCx_{\basisInd}(\vect x)$ and $\multiGPC_{\basisIndPar}(\vecpar)$ based on \reff{equ:indepBasis} and \reff{equ:nonGaussianBasis} independently. 
    \item[2.] Initialize the quadrature points for design variables $\{\vect x_l, v_l\}_{l=1}^{M_1}$ by (\ref{equ:independQuad}), and    quadrature points for stochastic parameters $\{\vecpar_l,u_l\}_{l=1}^{M_2}$ by the optimization problem \reff{equ:nonGaussianQuad}, respectively. Then co-optimize the quadrature rule to obtain $\{\vect x_k, \vecpar_k,w_k\}_{k=1}^{M}$ by (\ref{equ:co-optimQuad}). 
    \item[3.] Call the simulator to compute $f(\vect x_k, \vecpar_k)$, $y_i(\vect x_k, \vecpar_k)$ for all $i=1,\ldots,n$ and $k=1,\ldots,M$.
    \item[4.] Build the coefficients $h_{\basisInd,\basisIndPar}$ and  $c^i_{\basisInd,\basisIndPar}$ by equation \reff{equ:coeffs}.
    \item[5.] Set up the optimization problem \reff{equ:polyopt}, and then   solve it via a global polynomial optimization solver, e.g.,~\cite{henrion2009gloptipoly}.   
\end{enumerate}
\Output{The optimized design variable $\vect x^*$}
\end{algorithm}

\subsection{Basis Functions for Design and Uncertainty Variables}

For the mutually independent uniform-distributed  design variable  $\vect x$, their basis functions $\multiGPCx_{\basisInd}(\vect x)$ can be decoupled into the products of one-dimensional basis functions: 
\begin{equation}\label{equ:indepBasis}
    \multiGPCx_{\basisInd}(\vect x) = \phi_{\alpha_1}^1(x_1)\ldots \phi_{\alpha_{d_1}}^{d_1}(x_{d_1}).
\end{equation}
Here, $\phi_{\alpha_i}^i(x_i)$ is a Legendre polynomial~\cite{xiu2002wiener} and can be constructed by the three-term recurrence relation \cite{Walter:1982}. 

For the random vector $\vecpar$ describing non-Gaussian correlated process variations, we construct its basis functions $\multiGPC_{\basisIndPar}(\vecpar)$  by the Gram-Schmidt approach proposed in \cite{cui2018stochastic, cui2018stochastic2}. 
Specifically, we first reorder the monomials $\vecpar^{\basisIndPar}=\xi_1^{\beta_1}\ldots \xi_{d_2}^{\beta_{d_2}}$ in the  graded lexicographic order, and denote them as $\{p_j(\vecpar)\}_{j=1}^{N_{p}}$. Here, $N^{d_1}_{p}=\binom{d_2+p}{p}$ is the total number of basis functions for $\vecpar\in\mathbb{R}^{d_2}$ bounded by order $p$. 
Then we set $\multiGPC_1(\vecpar)   = 1$ and generate the orthonormal polynomials $\{\multiGPC_j(\vecpar)\}_{j=2}^{N_{p}}$ in the correlated parameter space recursively by
 \begin{align}\label{equ:nonGaussianBasis} \nonumber&\hat{\multiGPC}_j(\vecpar) = p_j(\vecpar)-\sum_{i=1}^{j-1} \mathbb{E}[ p_j(\vecpar)\multiGPC_i(\vecpar)] \multiGPC_i(\vecpar),
 \\
 &\multiGPC_j(\vecpar)  = \frac{\hat{\multiGPC}_j(\vecpar)}{\sqrt{\mathbb{E}[\hat{\multiGPC}^2_j(\vecpar)]}},\ j=2,\ldots,N_{p}.
 \end{align}
These basis functions $\{\multiGPC_j(\vecpar)\}_{j=1}^{N_{p}}$ can be reordered into  $\{\multiGPC_{\basisIndPar}(\vecpar)\}_{|\basisIndPar|=0}^{p}$.

\subsection{How to Compute the Coefficients?}
By a projection approach, the coefficient $ c^i_{\basisInd,\basisIndPar}$ for the basis function can be computed by 
\begin{equation}\label{equ:c}
    c^i_{\basisInd,\basisIndPar}= \mathbb{E}_{\vect x,\vecpar}[y_i(\vect x,\vecpar) \multiGPCx_{\basisInd}(\vect x)\multiGPC_{\basisIndPar}(\vecpar)].
\end{equation}
The above integration can be well computed given a suitable set of quadrature points $\{\vect x_k, \vecpar_k\}_{k=1}^{M}$ and weights $\{w_k\}_{k=1}^{M}$: 
\begin{align}\label{equ:coeffs}
c_{\basisInd,\basisIndPar}^i \approx\sum_{k=1}^{M}y_i(\vect x_k,\vecpar_k)\multiGPCx_{\basisInd}(\vect x_k)\multiGPC_{\basisIndPar}(\vecpar_k)w_k.
\end{align}
 We need to design a proper quadrature rule. The main challenge here is that $\vect x$ is an independent vector but $\vecpar$ describes  non-Gaussian correlated uncertainties.

 In this paper, we propose a three-stage optimization method to compute the quadrature rule:
 
\begin{itemize}
    \item Firstly, we compute the quadrature rule $\{\vect x_l, v_l\}_{l=1}^{M_1}$ for the independent design variable $\vect x$.
    
    \item Secondly, we employ the optimization approach proposed in~\cite{cui2018stochastic, cui2018stochastic2} to calculate the quadrature points and weights $\{\vecpar_l, u_l\}_{l=1}^{M_2}$ for the non-Gaussian correlated parameter $\vecpar$.
    
    \item Finally, we use their tensor products ($M_1M_2$ points) as an initialization and call the optimization approach proposed in~\cite{cui2018stochastic, cui2018stochastic2} for the coupled space of $\vect x$ and $\vecpar$ to compute $M\le M_1M_2$ joint quadrature points and weights $\{\vect x_k, \vecpar_k, w_k\}_{k=1}^{M}$. 
\end{itemize}  
The details are described below.



\subsubsection{Initial Quadrature Rule  for variables} One could employ the sparse grid approach \cite{sparse_grid:2000, Gerstner:1998} to compute  the quadrature samples and weights for the independent uniform random variables $\vect x\in\mathbb{R}^{d_1}$. However, the quadrature weights from a sparse grid method can be negative, and the number of quadrature points is not small enough. Therefore, after obtaining the sparse-grid quadrature rule, we  refine the quadrature rule by the least square optimization solver
\begin{equation}
\label{equ:independQuad}
\min_{\mat a\le\vect x_l\le \mat b, v_l\ge0}\  \sum_{j=1}^{N^{d_1}_{2p}}\left(\mathbb{E}[\multiGPCx_{j}(\vect x)]-\sum_{l=1}^{M_1}\multiGPCx_j(\vect x_l)v_l\right)^2.
\end{equation}
\ccf{Here, the expectations $\mathbb{E}[\multiGPCx_{j}(\vect x)]=\delta_{1j}$ are already known from the orthogonality of basis functions, and  $N^{d_1}_{2p}=\binom{d_1+2p}{2p}$.  This model is similar to that of \cite{cui2018stochastic, cui2018stochastic2}, which provides the quadrature points and weights to compute the numerical integral of all basis functions  upper bounded by order $2p$. 
If the optimized objective in  (\ref{equ:independQuad}) is small, the numerical integral of any functions in the $p$-th order polynomial space  will also be accurate. 
Further, the number of points $M_1$ can also be updated adaptively. 
The theoretical proofs for the number of quadrature points and the numerical approximation error  are provided in~\cite{cui2018stochastic2}.}
 
\subsubsection{Initial Quadrature Points for parameters} For the non-Gaussian correlated parameters $\vecpar$, we adopt the optimization-based quadrature rule in \cite{cui2018stochastic, cui2018stochastic2}. Specifically, we compute $M_2$ quadrature points $\vecpar_l$ and weights $w_l$ via  solving the following optimization problem
\begin{equation}
\label{equ:nonGaussianQuad}
\min_{\vecpar_l, u_l\ge0}\  \sum_{j=1}^{N^{d_2}_{2p}}\left(\mathbb{E}[\multiGPC_{j}(\vecpar)]-\sum_{l=1}^{M_2}\multiGPC_j(\vecpar_l)u_l\right)^2.
\end{equation}

\subsubsection{Optimized Joint Quadrature Points} 
The tensor product of the two sets of quadrature points $\{\vect x_l,v_l\}_{l=1}^{M_1}$ and  $\{\vecpar_l,u_l\}_{l=1}^{M_2}$ result in $M_1M_2$  simulation points in total, which may be still unaffordable for large-scale photonic design problems. We propose an optimization model to compute the joint quadrature rule for both the design variables $\vect x$ and the uncertain parameters $\vecpar$ to further reduce the simulation cost of building surrogate models:
\begin{align}
\label{equ:co-optimQuad}
\min_{\substack{\mat a\le\vect x_k\le \mat b\\ \vecpar_k, w_k\ge0}}\  \sum_{j_1=1}^{N^d_{2p}}\sum_{j_2=1}^{N^d_{2p}-j_1}&\left(\delta_{1j_1}\delta_{1j_2}-\sum_{k=1}^{M}\multiGPCx_{j_1}(\vect x_k)\multiGPC_{j_2}(\vecpar_k)w_k\right)^2.
\end{align}
Here $\delta_{1j_1}\delta_{1j_2}=1$ if $j_1=j_2=1$ and zero otherwise, and $d=d_1+d_2$. 
Our numerical experiments show that the total number of optimized quadrature points is $M$ is significantly smaller than $M_1M_2$. 
 
{\sl Remark:} Problem (\ref{equ:co-optimQuad}) is  a non-convex optimization  and is hard to optimize in general. The subproblems (\ref{equ:independQuad}) and (\ref{equ:nonGaussianQuad}) help to provide a good initial guess for the joint optimization.

We use the block coordinate descent optimization method  described in \cite{cui2018stochastic2} to solve all optimization subproblems (\ref{equ:independQuad}), (\ref{equ:nonGaussianQuad}), and (\ref{equ:co-optimQuad}).  
The following theorem ensures   high accuracy for our surrogate model considering the unavoidable numerical optimization error and function approximation error.

\begin{theorem}\cite{cui2018stochastic2}\label{thm:nonGaussian}
Assume that $\{\vect x_k, \vecpar_k,w_k\}_{k=1}^{M}$ are the numerical solution to \reff{equ:co-optimQuad}.
\begin{enumerate}
    \item Suppose \ccf{that} the objective function of \reff{equ:co-optimQuad} decays to zero. The required number of quadrature points is upper and lower bounded by
    \begin{equation}\label{equ:number_Qpoints}
        N^d_{p}=\frac{(d+p)!}{p! d!}\le M\le N^d_{2p}=\frac{(d+2p)!}{(2p)! d!};
    \end{equation}
    \item For any smooth and \zz{square-integrable} function $y(\vecpar)$, the approximation error of its $p$-th order stochastic approximation $\tilde y(\vecpar)$ satisfies 
    \begin{equation}\label{equ:sc_err}
\|y(\vect x, \vecpar)-\tilde y(\vect x, \vecpar)\|_2\le \alpha_1\delta_1 +\alpha_2\delta_2.
\end{equation}
Here, $\tilde{y}(\vect x, \vecpar)=\sum_{|\basisInd|+|\basisIndPar|=0}^{p}c_{\basisInd,\basisIndPar}\multiGPCx_{\basisInd}(\vect x)\multiGPC_{\basisIndPar}(\vecpar)$, $\delta_1$ is the $\ell_1$-norm of the objective function of \reff{equ:co-optimQuad} evaluated at its final numerical solution, $\delta_2$ is the distance of $y(\vect x,\vecpar)$ to the $p$-th order polynomial space, $\alpha_1=N_{p}LT$, $\alpha_2=1+N_{p}W$, $L=\max\|y(\vect x, \vecpar)\|_2$,  $T=\max_{j_1+j_2,l_1+l_2=1,\ldots,N_{2p}}\|\multiGPCx_{j_1}(\vect x)\multiGPC_{j_2}(\vecpar)\multiGPC_{l_1}(\vect x)\multiGPC_{l_2}(\vecpar)\|_2$, and $W=\sup \frac{|\mathbb I[y(\vecpar)]|}{\mathbb E[|y(\vecpar)|]}$ are constants\ccf{.}
\end{enumerate}

\end{theorem}

\ccf{ \textsl{Remark}: This subsection  focuses on the theory and implementation for building a surrogate model   for  low-dimensional problems.  For high-dimensional problems that are more costly in both computing the quadrature rule and difficult in reducing the number of samples, we may apply a high-dimensional solver such as the compressive  sensing~\cite{cui2019high} to build the surrogate model. Our framework  in Fig.~\ref{fig:flowchart} is still applicable.}

\begin{figure*}
  \begin{align}\label{equ:polyopt}
   \nonumber \min_{\vect x\in\ten X}\quad & \sum_{|\basisInd|=0}^{p}h_{\basisInd,\mat 0}\multiGPCx_{\basisInd}(\vect x)\\
  \text{s.t.}\quad &\kappa_{\epsilon_i}^2 \sum_{|\basisIndPar|=1}^{p}\left(\sum_{|\basisInd|=0}^{p-|\basisIndPar|}  c^i_{\basisInd,\basisIndPar}\multiGPCx_{\basisInd}(\vect x)\right)^2\le \left(u_i- \sum_{|\basisInd|=0}^{p}c^i_{\basisInd,\mat 0}\multiGPCx_{\basisInd}(\vect x)\right)^2,\quad  \sum_{|\basisInd|=0}^{p}c^i_{\basisInd,\mat 0}\multiGPCx_{\basisInd}(\vect x)\le u_i,\quad \forall\, i\in[n].
\end{align}
\end{figure*}

\subsection{The Proposed Polynomial Optimization}

With the formula for the mean value \reff{equ:mean} and the variance \reff{equ:var}, we obtain the following deterministic formula for the chance-constrained optimization:
  \begin{align}\label{equ:gPC_model}
   \nonumber \min_{\vect x\in\ten X}\ & \sum_{|\basisInd|=0}^{p}h_{\basisInd,\mat 0}\multiGPCx_{\basisInd}(\vect x)\\
  \nonumber  \text{s.t.}\ &\kappa_{\epsilon_i} \sqrt{ \sum_{|\basisIndPar|=1}^{p}\left(\sum_{|\basisInd|=0}^{p-|\basisIndPar|} c^i_{\basisInd,\basisIndPar}\multiGPCx_{\basisInd}(\vect x)\right)^2}\\
    &+\sum_{|\basisInd|=0}^{p}c^i_{\basisInd,\mat 0}\multiGPCx_{\basisInd}(\vect x)\le u_i,\ \forall\, i\in[n]. 
\end{align}
However, the constraints are non-smooth because of the square-root terms, and may not admit a gradient at some points~\cite{clarke1990optimization}. 
Instead, we use the equivalent smooth polynomial formula: 
\begin{align}
\kappa_{\epsilon_i}^2 \text{var}_{\vecpar}[y_i(\vect{x},\vecpar)]&\le (u_i- \mathbb{E}_{\vecpar}[y_i(\vect{x},\vecpar)])^2.
\end{align}
Consequently, \reff{equ:deterministic_model} can be reduced to a deterministic and smooth optimization problem of $\vect x$  in \reff{equ:polyopt}.

\begin{figure*}[t]
   \centering
    \includegraphics[width=6.6in]{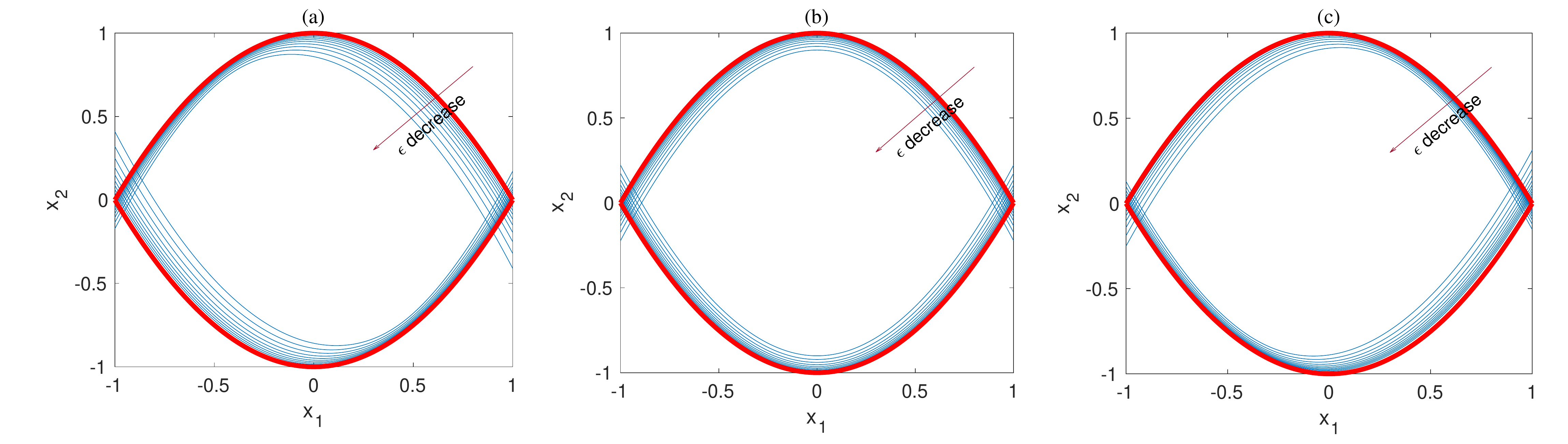}
    \caption{The feasible set of the synthetic example with risk tolerance levels  $\epsilon\in[10^{-2},10^{-0.1}]$ under different uncertainty distributions. (a): a positive-correlated non-Gaussian  distribution; (b): a Gaussian independent distribution; (c): a negative correlated non-Gaussian distribution. The domain between the red lines are the deterministic feasible set $x_1^2\pm x_2\le1$, and the blue lines demonstrate the  effects of chance constraints.}
    \label{fig:synthetic_feasibleset}
\end{figure*}

Noting that  both the objective function and the constraints  of  \reff{equ:polyopt} are polynomials, we can obtain the   optimal solution by using any polynomial solvers. In this paper, we  use the semi-definite relaxation based approaches \cite{waki2006sums, nie2013exact} because they can find the global optimal solution.

\section{Numerical Experiments}
In this section, we verify our proposed approach by a synthetic example and two photonic IC examples.  The p subproblem \reff{equ:polyopt} is solved by the global optimization solver GloptiPoly 3 \cite{henrion2009gloptipoly}. For a design variable $\vect x$, we generate $M$  parameters $\vecpar_j$ and approximate the yield  by
\begin{equation}\label{equ:yield_numerical}
    \text{yield}(\vect x) = \frac{\text{ is the number of } \vecpar_j \text{ such that} y_i(\vect x,\vecpar_j)\le u_i}{M}. 
\end{equation}  
We set all risk thresholds to $\epsilon$. For the synthetic example, we will compare our method with the deterministic formulation \reff{equ:P0_mean}. For the photonic IC examples, we will compare our method with the  Bayesian yield optimization method \cite{wang2017efficient}. We summarize the key idea of the Bayesian yield optimization in Appendix~\ref{append:BYO}. 
The MATLAB codes and a demo example can be downloaded online  \footnote{\url{https://web.ece.ucsb.edu/~zhengzhang/codes_dataFiles/ccyopt/}}.

\begin{figure*}[t]
    \centering
    \includegraphics[width=\textwidth]{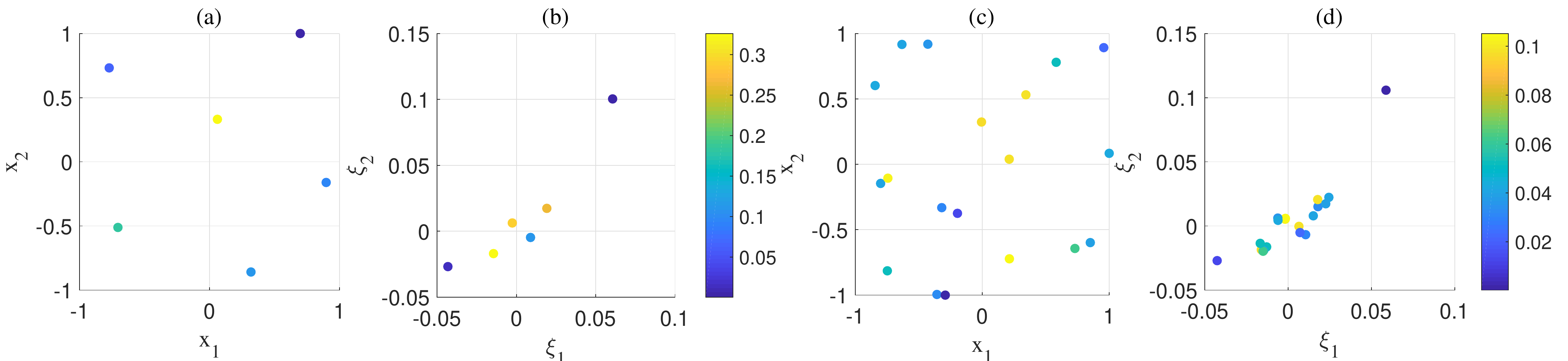}
    \caption{
    The quadrature points and weights in the synthetic experiment. (a) and (b): The initial 2-D quadrature points for the design variable $\vect x$ and uncertain parameters $\vecpar$ by solving (\ref{equ:independQuad}) and (\ref{equ:nonGaussianQuad}), respectively. (c) and (d): The optimized quadrature points for the joint 4-D space of $\vect x$ and $\vecpar$ by solving (\ref{equ:co-optimQuad}). Here we project the optimized 4-D quadrature points to the 2-D sub-space of $\vect x$ and $\vecpar$, respectively. The quadrature weights are shown in colors. }
    \label{fig:synthetic_Qpoints}
\end{figure*}

\subsection{Synthetic Example}\label{sec:synthetic}
Firstly, we consider a synthetic example with two design variables and two non-Gaussian correlated random parameters. The design variable $\vect x$ admits a uniform distribution $\mathcal{U}[-1,1]^2$ and the uncertain parameter $\vecpar$ follows a Gaussian mixture distribution. We define the yield criterion as $ (x_1+\xi_1)^2 \pm (x_2+\xi_2)\le 1$ and our goal is to maximize   $\mathbb{E}_{\vecpar}[3(x_1+\xi_1)+(x_2+\xi_2)]$. We formulate the yield into chance constraints and derive the following problem
\begin{align}\label{equ:synthex}
  \nonumber  \max_{\vect x}\quad& \mathbb{E}_{\vecpar}[3(x_1+\xi_1)-(x_2+\xi_2)] \\
  \nonumber   \text{s.t.}\quad &\text{Prob}_{\vecpar}\left( (x_1+\xi_1)^2 -(x_2+\xi_2)\le 1\right)\ge 1-\epsilon,\\
     &\text{Prob}_{\vecpar}\left( (x_1+\xi_1)^2 +(x_2+\xi_2)\le 1\right)\ge 1-\epsilon.
\end{align}

To  illustrate  the effects of different parameter distributions, we study three probability density functions: the independent distribution $\ten N(\mat 0, 10^{-4}\mat I)$, the    non-Gaussian positive correlations $\frac12\ten N(\mat {0.01},10^{-4}\boldsymbol{\Sigma})+\frac12\ten N(-\mat {0.01},10^{-4}\boldsymbol{\Sigma})$ with $\boldsymbol{\Sigma}= \left(\begin{array}{cc}
1 & 0.75\\0.75& 1\end{array}\right)$, and the non-Gaussian negative correlations $\frac12\ten N([0.01,-0.01]^T,10^{-4}\boldsymbol{\Sigma})+\frac12\ten N([-0.01,0.01]^T,10^{-4}\boldsymbol{\Sigma})$ with $\boldsymbol{\Sigma}= \left(\begin{array}{cc}
1 & -0.75\\-0.75& 1\end{array}\right)$. 
The feasible sets under three probability density distributions are shown in  Fig.~\ref{fig:synthetic_feasibleset}.  The comparison clearly shows that the effects of different uncertainties. For all three density functions, the feasible regions are reduced when the risk level $\epsilon$ decreases.

 \begin{table}[t]
  
     \centering
          \caption{The optimal solution for the synthetic experiment under difference risk threshold $\epsilon$.} 
    \begin{tabular}{|c|cc|c|c|}
    \hline
   Algorithm &\multicolumn{2}{c|}{$\vect x^*$} & Objective & Yield (\%) \\  \thickhline
   Proposed ($\epsilon=$ 0.01) &  0.8630 &  -0.1172  &  2.4717  &100\\
     Proposed ($\epsilon=$  0.05) &    0.9379  & -0.0522   & 2.7616 & 100\\
      Proposed ($\epsilon=$ 0.10)&   0.9587 &  -0.0402 &   2.8360  & 99.42\\
      Proposed ($\epsilon=$ 0.15) &  0.9689  & -0.0351   & 2.8717  & 93.84\\
       Proposed ($\epsilon=$0.20)&    0.9751 &  -0.0293   & 2.8959 &  87.49\\\hline 
 \reff{equ:P0_mean} & \ccf{0.9999}     &    0   & \ccf{2.9997}  &  \ccf{41.66}\\
     \hline 
    \end{tabular}
    \label{tab:synthetic}
 \end{table}

Next we take the non-Gaussian positive correlated distribution as an example  to compute the optimal solution of \reff{equ:synthex}. We first build the surrogate models for both the objective and constraints by the second-order polynomial basis functions. The optimized quadrature points $\{\vect x_l, v_l\}_{l=1}^6$ for the design variables by \reff{equ:independQuad} and $\{\vecpar_l, u_l\}_{l=1}^6$ for the random parameter by (\ref{equ:nonGaussianQuad}) are shown in Fig.~\ref{fig:synthetic_Qpoints}~(a) and (b), respectively. 
Directly tensorizing the two sets of quadrature points generates $36$ samples. We further solve  (\ref{equ:co-optimQuad}) to reduce them to $M=19$ optimized samples and weights. According to  Theorem~\ref{thm:nonGaussian}, the number of quadrature samples for $d=4, p=2$ should be in the range $[15,70]$. Our optimization algorithm obtains $M=19$, which is close to the theoretical lower bound.

\begin{figure}[t]
    \centering
    \includegraphics[width=0.4\textwidth]{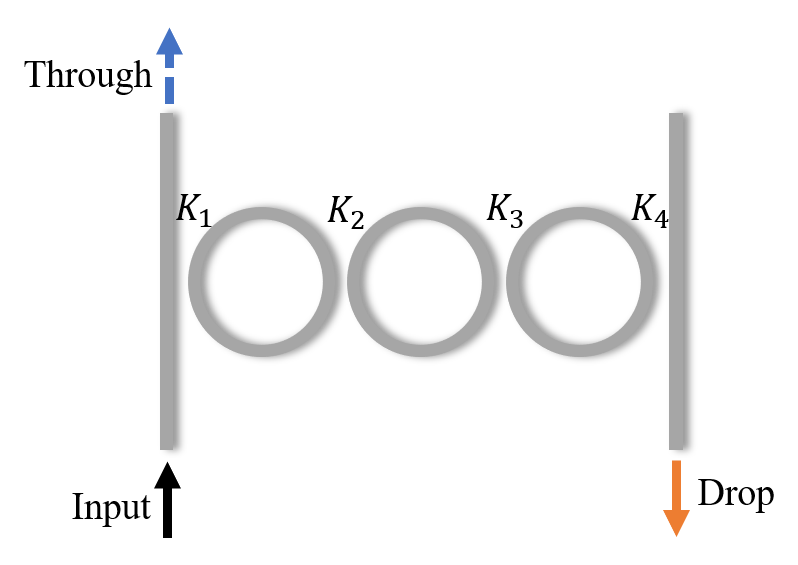}
    \caption{An optical add-drop filter with three microrings coupled in series. }
 
    \label{fig:filter_scheme}
\end{figure}

We further show the results for  different risk tolerance levels $\epsilon$ in Table~\ref{tab:synthetic}.  A smaller $\epsilon$ results in a smaller feasible domain (as shown in Fig.~\ref{fig:synthetic_feasibleset}), and generates a higher yield but a smaller objective value.  \ccf{In practice, $\epsilon$ can be chosen case-by-case based on the trade-off between the performance and yield requirements.}  Compared with the solution $\tilde{\vect x}=[0.9999,0]^T$ from solving \reff{equ:P0_mean}, our method can achieve a significantly higher yield: our optimized yield is above $87\%$ while solving \reff{equ:P0_mean} only leads to a yield of $41.66\%$. 


\subsection{Microring \ccf{Add-drop} Filter}
\label{sec:experiment_band_pass}

\begin{figure*}[t]
    \centering
    \includegraphics[width=\textwidth]{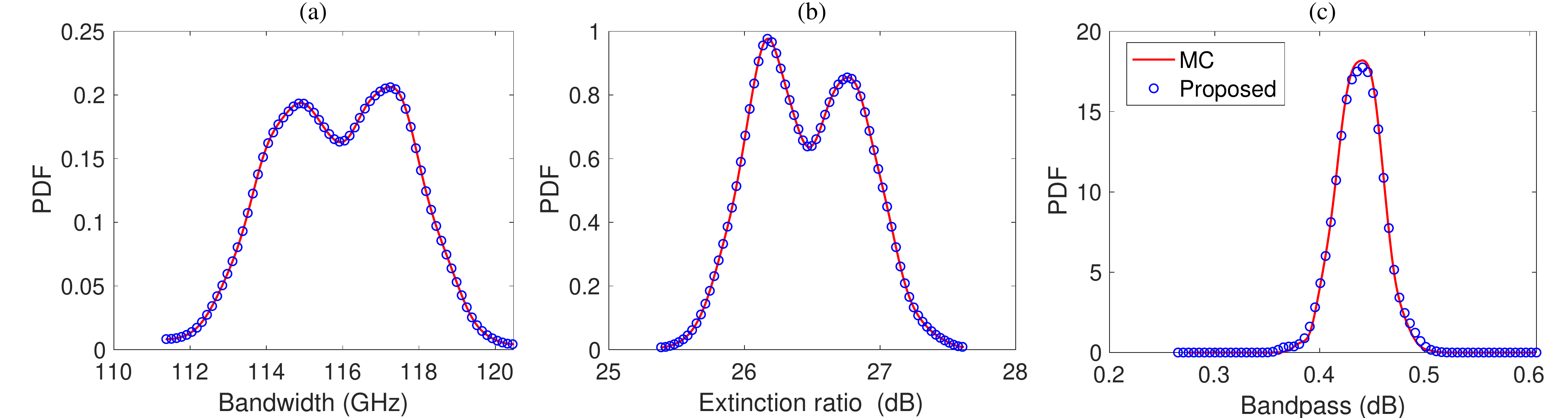}
    \caption{The probability density functions (PDF) of the bandwidth, extinction ratio and roughness for the microring add-drop filter at the optimal solution $\vect x^*=[0.5582,  0.4208, 0.3000, 0.6000]$ by our proposed optimization with  $\epsilon=0.05$. Our surrogate model uses only 64 simulations, and Monte Carlo (MC) uses $10^3$ simulations.}
    \label{fig:gPC_app_err_ring}
\end{figure*}

\begin{figure*}[t]
    \centering
        \includegraphics[width=\textwidth]{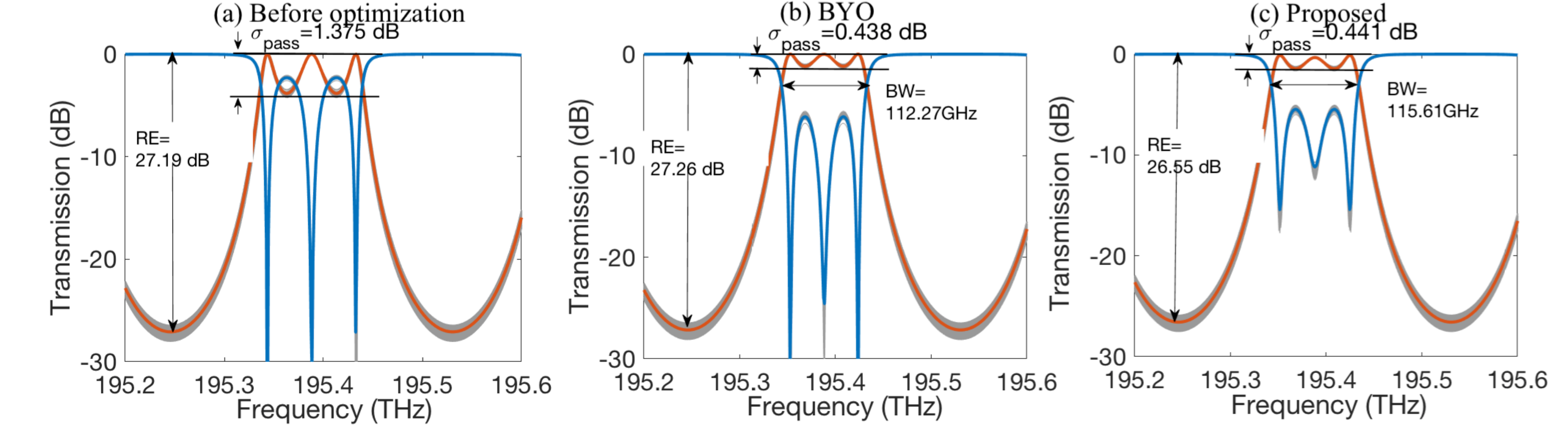}
    \caption{The transmission curves of the microring add-drop filter at different design choices. The grey lines show the uncertainties caused by the process variations. The orange and blue curves show the mean transmission rates at the drop port and the through port, respectively. Here $\text{RE}$, $\text{BW}$ and $\sigma_{pass}$ denote the mean values of extinction ratio, bandwidth and roughness, respectively.
    (a) The transmission at $\vect x^0=[0.45,0.45,0.45,0.45]$  without any optimization. It doesn't have a clear passband because $\sigma_{pass}$ is too large. (b) The results after the Bayesian yield optimization; (c) The results obtained from our chance-constrained optimization with $\epsilon=0.05$.  
    }
 
    \label{fig:filter}
\end{figure*}

We continue to consider the design of an optical add-drop filter consisting of three identical silicon microrings coupled in series, as shown in Fig. \ref{fig:filter_scheme}. In designing such a broadband optical filter, the coupling coefficients play an important role in determining  key performance metrics, such as the bandwidth and extinction ratio\cite{orta1995synthesis, pintus2013analysis}. A broad and flat passband with a high extinction ratio can be achieved by optimizing the coupling strengths between the microrings \cite{orta1995synthesis}. In this example, we employ silicon as the waveguide material and assume the effective refractive index to be $n_{\text{eff}}=2.44$ and the effective group index to be $n_\text{g}=4.19$ near the wavelength of 1.55 $\mu m$. The design variables are  the coupling coefficients $\vect{x}=[K_1,K_2,K_3,K_4]$ that are to be optimized within the interval of $[0.3,0.6]$.
The random variables are set as  small deviations of the coupling coefficients. We assume that $\vect{\xi}$ follows a non-Gaussian correlated distribution \ccf{
\begin{equation}
   p(\vecpar) =  \frac12 \ten N(\vect \mu_1,\vect \Sigma) + \frac12 \ten N(\vect \mu_2,\vect \Sigma),
\end{equation}
where $\vect \mu_1=-\vect \mu_2=0.006[1,1,1,1]^T$, 
and the variance is defined as 
$\vect \Sigma = 0.006^2\left[
    \begin{array}{cccc}
       1  & 0.4 & 0.1 & 0.4 \\
        0.4 & 1 & 0.4 & 0.1 \\
        0.1 & 0.4 & 1 & 0.4\\
        0.4 & 0.1 & 0.4 & 1
    \end{array}
    \right].$}

We mainly focus on three metrics of the microring filter: the 3dB bandwidth ($\text{BW}$, in GHz), the extinction ratio ($\text{RE}$, in dB) of the transmission at the drop port, and the roughness ($\sigma_{pass}$, in dB) of the passband that takes a standard deviation of the passband. The yield-aware optimization problem of the microring filter design can be formulated as:

\begin{align}\label{equ:opt_ring}
  \nonumber  \max_{\vect x\in\ten X}\quad& \mathbb{E}_{\vecpar}[\text{BW}(\boldsymbol{x},\boldsymbol{\xi})]\\
  \nonumber  \text{s.t.}\quad&            \text{Prob}_{\boldsymbol{\xi}}(\text{RE}(\boldsymbol{x},\boldsymbol{\xi})\ge \text{RE}_0)\ge 1-\epsilon,\\
             &\text{Prob}_{\boldsymbol{\xi}}(\sigma_{pass}(\boldsymbol{x},\boldsymbol{\xi})\le \sigma_0)\ge 1-\epsilon,
\end{align}
where the yield is defined via some chance constraints on the extinction ratio and the roughness of the passband. In our simulation, the threshold extinction ratio (RE$_0$) and the roughness of the passband ($\sigma_0$) are 25dB and 0.5dB, respectively.

\begin{table}[t]
    \centering
    \caption{Optimization results for the microring add-drop filter.}
    \begin{tabular}{|c|c|c|c|}
    \hline
 Algorithm & Simulations& $\mathbb{E}_{\vecpar}[\text{BW}]$ (GHz) & Yield    (\%)\\
  \thickhline

Proposed ($\epsilon=0.03$)&64   &  113.4 & 100 \\
Proposed ($\epsilon=0.05$)&64   &  115.6 & 99.8 \\
Proposed ($\epsilon=0.07$)&64   &  117.2 & 99.5 \\
Proposed ($\epsilon=0.10$)&64   &  118.4 & 98.1 \\
\hline
BYO~\cite{wang2017efficient}      & 2020 & 112.3 & 99.8 \\ \hline
    \end{tabular}
    \label{tab:comp_ring}
 \end{table}
 
 We first build the second-order polynomial surrogate model by our proposed Algorithm~\ref{alg:framework}. We only need 17 initial quadrature points for the variable $\vect x$ by solving \reff{equ:independQuad}, 16 quadrature points for the parameters $\vecpar$ by solving \reff{equ:nonGaussianQuad}, and 64 quadrature points for the joint optimization of $\vect x$ and $\vecpar$ by solving \reff{equ:co-optimQuad}.   Fig.~\ref{fig:gPC_app_err_ring} shows that our surrogate model can well approximate the probabilistic distributions of the performance metrics with the comparison of $10^3$ Monte Carlo simulations, although our method only needs 64 simulation samples for this example.

\begin{figure}[t]
    \centering
    \includegraphics[width=0.5\textwidth]{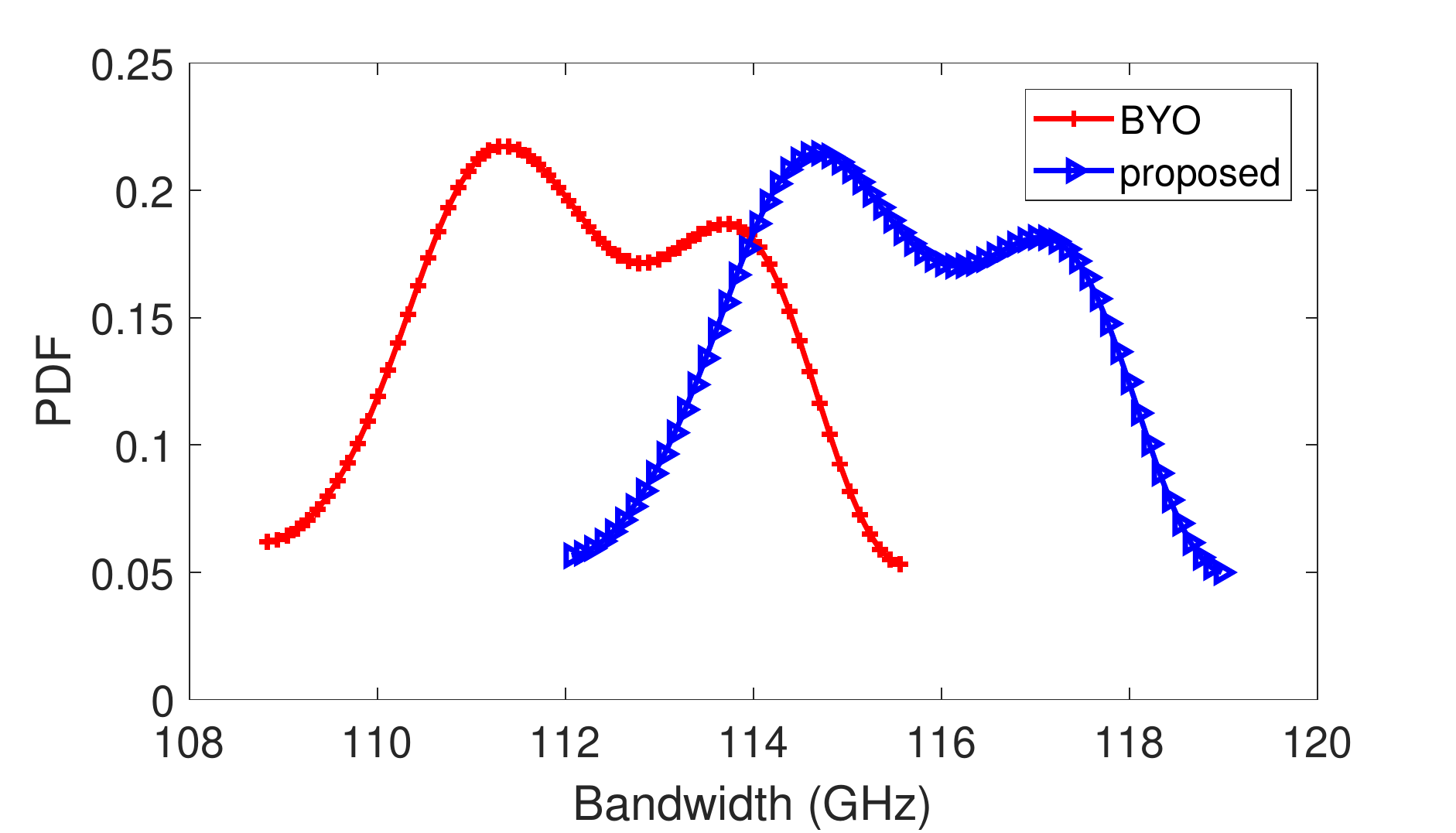}
     \caption{The optimized bandwidth probability density distribution of the microring filter. Our chance-constrained optimization obtain an expected value of 115.6 GHz while the Bayesian yield optimization (BYO) only produces an expected value of 112.3 GHz. 
    }
    \label{fig:BW_filter_MC}
\end{figure}

We  summarize the results of our proposed method with different choices of $\epsilon$ and the results obtained by the Bayesian yield optimization (BYO) in Table~\ref{tab:comp_ring}. 
It shows that when risk tolerance level $\epsilon$ decreases, our proposed method can achieve higher yield and lower bandwidth. This is corresponding to our theory that   a lower risk level $\epsilon$   results in a smaller   feasible region. Our proposed method can always achieve a large bandwidth because it computes the global optimal solution of the polynomial optimization problem. 
\ccf{When $\epsilon=0.05$, we get a bandwidth  $\mathbb{E}_{\vecpar}[\text{BW}]=115.6$ GHz with   $99.8\%$ yield at the   optimal solution $\vect{x}^*=[0.5582,0.4208,0.3000,0.6000]$, while BYO takes 2020 simulations to achieve the result of $\mathbb{E}_{\vecpar}[\text{BW}]=112.3$ GHz with the yield 99.8\%.} 
 Fig.~\ref{fig:filter} compares the \ccf{frequency response} before and after the yield-aware optimization. Both our proposed method and BYO can achieve a higher bandwidth with a smoother passband  compared to the design before optimization.  
 In Fig. \ref{fig:BW_filter_MC}, we  plot the probability density of the bandwidth at the optimal design by our chance-constrained optimization with $\epsilon=0.05$ and by the BYO, respectively. It clearly shows that our proposed method can increase the bandwidth while achieving the same yield. 


\begin{figure}[t]
    \centering
    \includegraphics[width=0.5\textwidth]{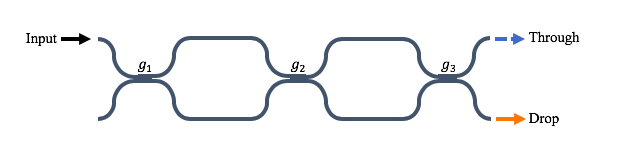}  
    \caption{The schematic of a third-order Mach-Zehnder Interferometer.}
    \label{fig:MZ_scheme}
\end{figure}

\subsection{Mach-Zehnder Interferometer}
\label{sec:experiment_MZ}

We apply the same framework to optimize a third-order Mach-Zehnder interferometer (MZI) which consists of three port coupling and two  arms, as shown in Fig. \ref{fig:MZ_scheme}.  The coupling coefficients between the MZ arms play the most important role in the design. 
\ccf{The relationship between the coupling coefficient $\kappa$ and the gap $g$ (nm) is } 
\begin{equation}
    \kappa = \exp(-\frac{g}{260}).
\end{equation}
\ccf{In this experiment, the design variables $\vect x=[g_1,g_2,g_3]$   are  optimized  in the interval of [100 nm, 300 nm]$^3$.}
The random variable $\vect{\xi}$  follows the  Gaussian mixture  distribution 
\ccf{\begin{equation}
   p(\vecpar) =  \frac12 \ten N(\vect \mu_1,\vect \Sigma) + \frac12 \ten N(\vect \mu_2,\vect \Sigma),
\end{equation}
where $\vect \mu_1=-\vect \mu_2=[1,1,1]^T$,  
and $\vect \Sigma =\left[
    \begin{array}{ccc}
       1  & 0.4 & 0.1  \\
        0.4 & 1 & 0.4  \\
        0.1 & 0.4 & 1 \\
    \end{array}
    \right].$}
We consider three performance metrics of the MZI: the 3dB bandwidth ($\text{BW}$, in GHz), the crosstalk ($\text{XT}$, in dB), and the attenuation ($\alpha$, in dB) of the peak transmission. The yield is defined through the crosstalk and the attenuation. The yield-aware optimization is formulated as:
\begin{align} \label{equ:opt_mz}
  \nonumber  \max_{\vect x}\quad& \mathbb{E}_{\vecpar}[\text{BW}(\vect{x},\boldsymbol{\xi})]\\
  \nonumber  \text{s.t.}\quad&  \text{Prob}_{\boldsymbol{\xi}}(\text{XT}(\boldsymbol{K},\boldsymbol{\xi})\le \text{XT}_0)\ge 1-\epsilon,\\
             &\text{Prob}_{\boldsymbol{\xi}}(\alpha(\boldsymbol{x},\boldsymbol{\xi})\le \alpha_0)\ge 1-\epsilon,
\end{align}
where the yield risk level is  $\epsilon$. The threshold crosstalk ($\text{XT}_0$) and attenuation ($\alpha_0$) are -4 dB and  2 dB, respectively.

\begin{figure*}[t]
    \centering
    \includegraphics[width=0.9\textwidth]{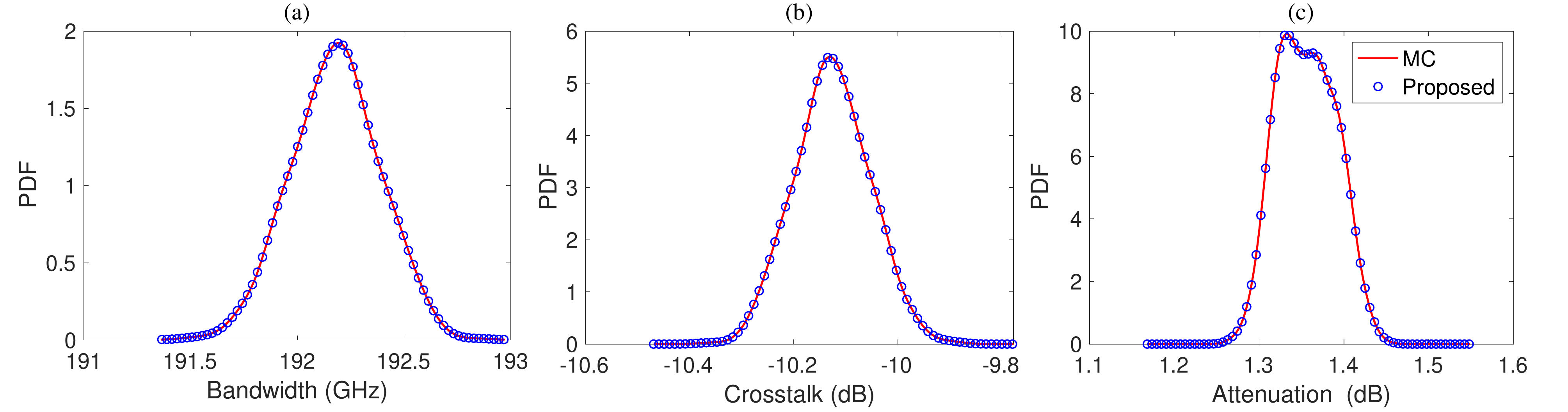}
    \caption{The probability density functions (PDF)  for the bandwidth, crosstalk, and attenuation of the MZI at our optimized design parameters $\vect x^*=[0.300, 0.5036,0.300]$. Our surrogate model uses only 36 simulations and Monte Carlo (MC) uses $1000$ simulations. }
    \label{fig:gPC_app_err_MZ}
\end{figure*}

 We first build three second-order polynomial surrogate models for $\text{BW}$, $\text{XT}$ and $\alpha$ by our proposed Algorithm~\ref{alg:framework}. We generate 11 initial quadrature points for the design variable $\vect x$, 10 initial quadrature points for the uncertainty parameter $\vecpar$. Then we apply the tensor product of those 110 points to problem (\ref{equ:co-optimQuad})  and eventually get $36$ quadrature points for the joint space after   co-optimization.  Fig.~\ref{fig:gPC_app_err_MZ} shows that our surrogate models constructed with 36 quadrature points can well approximate the density functions of all three performance metrics compared with Monte Carlo with $10^3$ samples.

\begin{figure*}[t]
    \centering
    \includegraphics[width=0.9\textwidth]{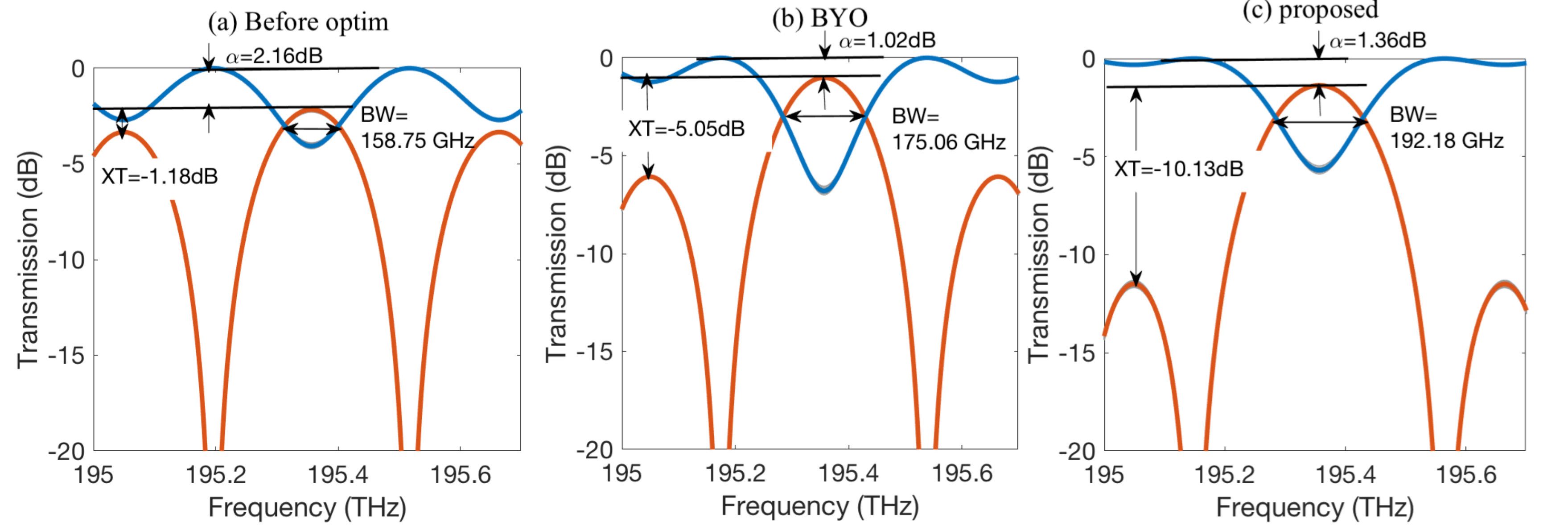}
    \caption{The transmission curves of the MZI. The grey lines show the performance uncertainties. The orange and blue curves show the transmission rates at the drop and through ports, respectively. The mean values of the bandwidth, crosstalk and attenuation are denoted as $\text{BW}$, $\text{XT}$ and $\alpha$, respectively.
    (a) The initial design  $\boldsymbol{x}^{0}=[150,150,150]$; (b) Design after Bayesian yield optimization; (c) Design with the proposed chance-constrained yield optimization.}
    \label{fig:MZ_filter}
\end{figure*}

\begin{figure}[t]
    \centering
    \includegraphics[width=0.48\textwidth]{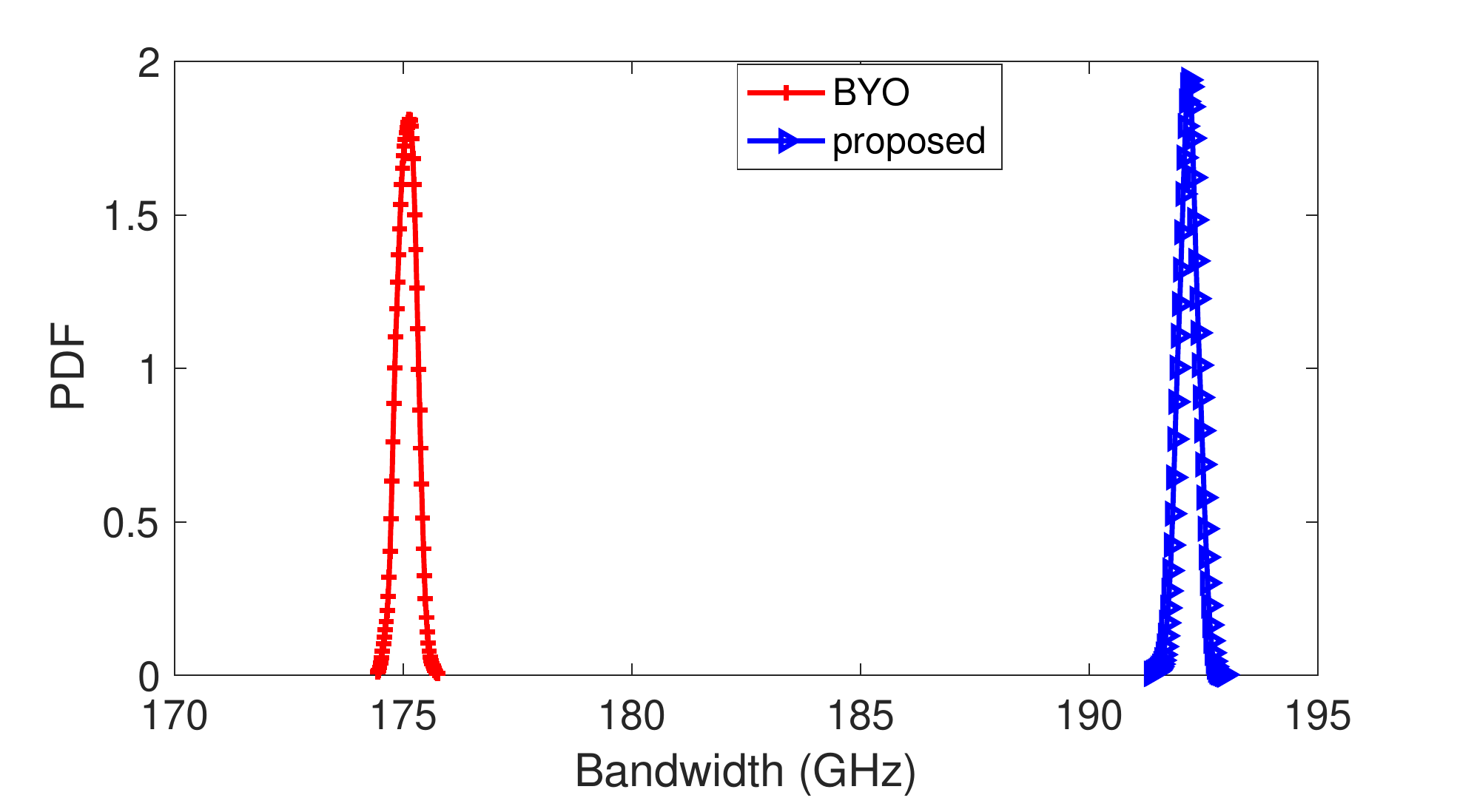}
    \caption{The optimized bandwidth of the MZI by the Bayesian yield optimization and our proposed method, respectively. The expectation bandwidth of the Bayesian yield optimization is 175.4 GHz while our proposed method with $\epsilon=0.05$ can get 186.4 GHz.}
    \label{fig:MZ_MC}
\end{figure}

We also compare our proposed method and BYO in Table~\ref{tab:comp_MZ}. Similar to the result in Table~\ref{tab:comp_ring}, a lower risk tolerance results in higher yield and a lower  expected value of bandwidth. Our method requires $56\times$ fewer simulation points than BYO, which is a great advantage for design cases with the time-consuming simulations.
For $\epsilon=0.05$, the optimized nominal design is 
\ccf{$\vect x^*=[300, 111.2, 300]$} and its expected bandwidth is 192.2 GHz.
In Fig. \ref{fig:MZ_filter}, we  compare the \ccf{frequency response} before and after the yield-aware optimization. Our proposed method can have a higher bandwidth and a smaller crosstalk  compared to Bayesian yield optimization and the initial design. Fig. \ref{fig:MZ_MC} further shows the probability density of the optimized bandwidth by our chance-constrained optimization and the Bayesian yield optimization, respectively. It clearly shows that our proposed method produces higher bandwidth. 

\begin{table}[t]
    \centering
    \caption{Optimization result for the MZI.}
    \begin{tabular}{|c|c|c|c|}
\hline
  Algorithm  & Simulations & $\mathbb{E}_{\vecpar}[\text{BW}]$ (GHz) & Yield (\%) \\
  \thickhline

   Proposed ($\epsilon=0.03$) &  36  &  188.8  &  100 \\
   Proposed ($\epsilon=0.05$) &  36  &  192.2 &  100 \\
   Proposed ($\epsilon=0.07$) &  36  &  194.5  &  100 \\
  Proposed ($\epsilon=0.10$) &  36  &  195.0  & 87.7 \\
        \hline
      BYO~\cite{wang2017efficient}      &  2020 &  175.0  &  100 \\\hline
   
    \end{tabular}
    \label{tab:comp_MZ}
 \end{table}

\section{Conclusions and Remarks}
This paper has presented a data-efficient framework for the yield-aware optimization of \ccf{photonic ICs} under non-Gaussian correlated process variations. We have proposed to reformulate the stochastic chance-constrained optimization into a deterministic polynomial optimization problem. Our framework only requires simulation at a small number of important points and admits a surrogate model for yield-aware optimization. In the experiments by the microring filter and the Mach Zehnder filter, we have demonstrated that our optimization scheme can give   high yield and high bandwidth. Compared with Bayesian yield optimization, our method has consumed much fewer simulation samples and produced better design performance while achieving the same yield. 

This work should be regarded as a presentation of  preliminary results in this direction. Many problems are worth further investigation in the future, for instance:
\begin{itemize}

\item { Non-Smoothness.} Similar to generalized polynomial chaos~\cite{xiu2002wiener}, the surrogate modeling techniques in~\cite{cui2018stochastic,cui2018stochastic2} require the stochastic functions to be smooth. However,  performance metrics of a photonic IC may be non-smooth with respect to the design variables and process variations. How to handle non-smoothness in this optimization framework is a critical issue.
    
    \item {High Dimensionality.} Large-scale photonic ICs may have a huge number of design variables and process variation parameters. This brings new challenges to the surrogate modeling and the resulting polynomial optimization in our framework.
    
\end{itemize}

\appendices

\section{Detailed Derivation of Equation \reff{equ:chance2pop}}
\label{app:chance_proof}


Suppose $u> \mathbb{E}_{\vecpar}[y(\vect{x},\vecpar)]$. We show  that  the following deterministic constraint 
\begin{equation*} 
    \mathbb{E}_{\vecpar}[y(\vect{x},\vecpar)]+\sqrt{(1-\epsilon)/\epsilon} \sqrt{\text{var}_{\vecpar}[y(\vect{x},\vecpar)]}\le u
\end{equation*}
is a sufficient but not necessary condition for the   chance constraint: 
\begin{equation*} 
    \text{Prob}_{\vecpar}(y(\vect x,\vecpar)\le  u)\ge 1-\epsilon.
\end{equation*}
 In other words, we want to show that each feasible point of \reff{equ:determin_chance} is a feasible point of the chance constraint \reff{equ:yield2chance}.

 Denote the random variable as $\vect X=y(\vect{x},\vecpar)$.   Cantelli's inequality \cite{cantelli1929sui} states that for any random variable $\vect X$ with a mean value $\mathbb{E}[\vect X]=\mathbb E_{\vecpar}[y(\vect x,\vecpar)]$ and variance $\sigma^2=\text{var}_{\vecpar}[y(\vect x,\vecpar)]$, it holds that the probability of a single tail can be bounded as follows:
  \begin{equation}
      \text{Prob}(\vect X- \mathbb{E}[\vect X]\le \lambda ) \ge 1-\frac{\sigma^2}{\sigma^2+\lambda^2} \ \text{ if } \lambda>0. 
  \end{equation}
Therefore, for any constant $u\ge \mathbb E[\vect x]$ we have 
  \begin{align*} 
      \text{Prob}(\vect X\le u )&=\text{Prob}(\vect X-\mathbb E[\vect X]\le u-\mathbb E[\vect X] ) \\
      &\ge 1-\frac{\sigma^2}{\sigma^2+( u-\mathbb E[\vect X])^2}. 
  \end{align*} 
For any $\epsilon$, a sufficient condition for $\text{Prob}(\vect X\le u )\ge 1-\epsilon$ is $1-\frac{\sigma^2}{\sigma^2+( u-\mathbb E[\vect x])^2}\ge 1-\epsilon$, i.e., 
\begin{align}
    \mathbb{E}[\vect   X]+\sqrt{(1-\epsilon)/\epsilon}\sigma\le u.
\end{align}
Substituting $\vect X=y(\vect x,\vecpar)$ into the above equation we get \reff{equ:chance2pop}. The proof is completed.

\section{Detailed derivation of  equations (\ref{equ:mean}) and (\ref{equ:var})} 
\ccf{
Suppose that the smooth function $y(\vect x,\vecpar)$ \zz{is already} represented by a linear combination of \zz{some} basis functions,}
\begin{equation}\label{equ:surrogate_append}
    y(\vect x,\vecpar)=\sum_{|\basisInd|+|\basisIndPar|=0}^{p} c_{\basisInd,\basisIndPar}\multiGPCx_{\basisInd}(\vect x)\multiGPC_{\basisIndPar}(\vecpar),
\end{equation}
where  $\mathbb E[\multiGPC_{\basisIndPar}(\vecpar)\multiGPC_{\boldsymbol \gamma}(\vecpar)]=\delta_{\basisIndPar,\boldsymbol{\gamma}}$. The mean value of $y(\vect x,\vecpar)$ is 
\begin{align*} 
    \mathbb{E}_{\vecpar}[y(\vect x,\vecpar)] &= \sum_{|\basisInd|=0}^{p}\sum_{|\basisIndPar|=0}^{p-|\basisInd|}c_{\basisInd,\basisIndPar}\multiGPCx_{\basisInd}(\vect x)\mathbb E[\multiGPC_{\basisIndPar}(\vecpar)]\\
    &= \sum_{|\basisInd|=0}^{p}c_{\basisInd,0}\multiGPCx_{\basisInd}(\vect x),
\end{align*}
where the last equality is due to $\multiGPC_{0}(\vecpar)=1$ and $\mathbb E[\multiGPC_{\basisIndPar}(\vecpar)]=\mathbb E[\multiGPC_{\basisIndPar}(\vecpar)\multiGPC_{0}(\vecpar)]=0$, $\forall \basisIndPar\neq 0$.
The variance is  
\begin{align*} 
    \text{var}_{\vecpar}[y(\vect x,\vecpar)]
    &=\mathbb E_{\vecpar}[(y(\vect x,\vecpar)-\mathbb E_{\vecpar}[y(\vect x,\vecpar)])^2]\\
    &= \mathbb E_{\vecpar}\left[ \left(\sum_{|\basisIndPar|=1}^{p}
    \left( 
    \sum_{|\basisInd|=0}^{p-|\basisIndPar|}c_{\basisInd,\basisIndPar}\multiGPCx_{\basisInd}(\vect x)
    \right)\multiGPC_{\basisIndPar}(\vecpar)\right)^2\right]\\
    &= \sum_{|\basisIndPar|=1}^{p}\left(\sum_{|\basisInd|=0}^{p-|\basisIndPar|} c_{\basisInd,\basisIndPar}\multiGPCx_{\basisInd}(\vect x)\right)^2, 
\end{align*}
where the last equality is due to the basis functions $\{ \multiGPC_{\basisIndPar}(\vecpar)\}$ are orthogonal in the stochastic parameter space.

  \section{Bayesian Yield Optimization (BYO)}
  \label{append:BYO}
Bayesian yield optimization (BYO) is a state-of-the-art tool for the yield optimization of electronic devices and circuits~\cite{wang2017efficient}. 
This method approximates and optimizes the posterior distribution of design variable under the condition of ``pass'' events 
\ccf{$$S=\{(\vect x, \vecpar): (\vect x, \vecpar) \text{ satisfies all performance constraints}\}.$$} 
With the Bayes' theorem, it holds that  $\text{Prob}(S|\vect{x})=\frac{\text{Prob}(S)}{\text{Prob}(\vect{x})}\text{Prob}(\vect{x}|S)$. 
In our problem setting,  $\text{Prob}(\vect{x})$ is a constant because we assume that $\vect x$ follows a uniform distribution and $\text{Prob}(S)$ should also be a constant without the dependence on the variable $\vect x$. Therefore,   $\text{Prob}(S|\vect{x})\propto \text{Prob}(\vect{x}|S)$ and  the original yield optimization problem \reff{equ:max_yield} is equivalent to 
\begin{equation}\label{equ:max_yield4}
    \boldsymbol{x}^{BYO}= \underset{\vect x\in \mathcal{X}}{\text{argmax}} \ \ \text{Prob}(\boldsymbol{x}|S). 
\end{equation}
The paper \cite{wang2017efficient} proposed an expectation-maximization framework to solve problem (\ref{equ:max_yield4}). 
At the $t$-th iteration,  the expectation step approximates the probability by the kernel density estimation. 
Specifically,  we generate $N=100$ samples $(\vect x_i, \vecpar_i)$ randomly  and  call a simulator to compute the quantity of interests at those samples. Then   choose $M\le N$ ``pass" samples to perform the kernel density estimation 
\begin{equation*}
    \text{Prob}(\vect{x}|S)\approx\frac{1}{M}\sum_{i=1}^{M}\frac{1}{\sqrt{2\pi}h}\exp{(-\frac{1}{2h}(\vect x-\vect \mu_i)^T(\vect x-\vect \mu_i))}, 
\end{equation*}
where $\{\vect \mu_i\}_{i=1}^M\in S$ are  design samples that  satisfies the performance constraints and $h=0.3$  is a bandwidth parameter. 
Afterward, the maximization step returns an updated design variable $\vect x^{BYO,t}$. 
We will call the simulator again at this design variable to record its objective value and ``pass" status. 
We terminate the   algorithm if the maximal iteration number $20$ is reached, or the residue of two consecutive iterations is below $10^{-6}$. 
After the whole optimization process, we return the design variable that can pass the yield constraints  with the best objective value
\begin{equation*}
   \vect x^{BYO}=\arg\max_{\vect x \in \vect x^{BYO,t}} \mathbb{E}_{\vecpar}[f(\vect x,\vecpar)]\text{ s.t. }  \text{pass}(\vect x)=1.
\end{equation*}

\section*{Acknowledgment}
The authors would like to thank the anonymous reviewers for their detailed comments. 
We also appreciate Paolo Pintus for his helpful discussions on the benchmarks.

 \bibliographystyle{IEEEtran}
\bibliography{ref.bib}

 \begin{IEEEbiography}
  [{\includegraphics[width=1in,height=1.25in,clip,keepaspectratio]{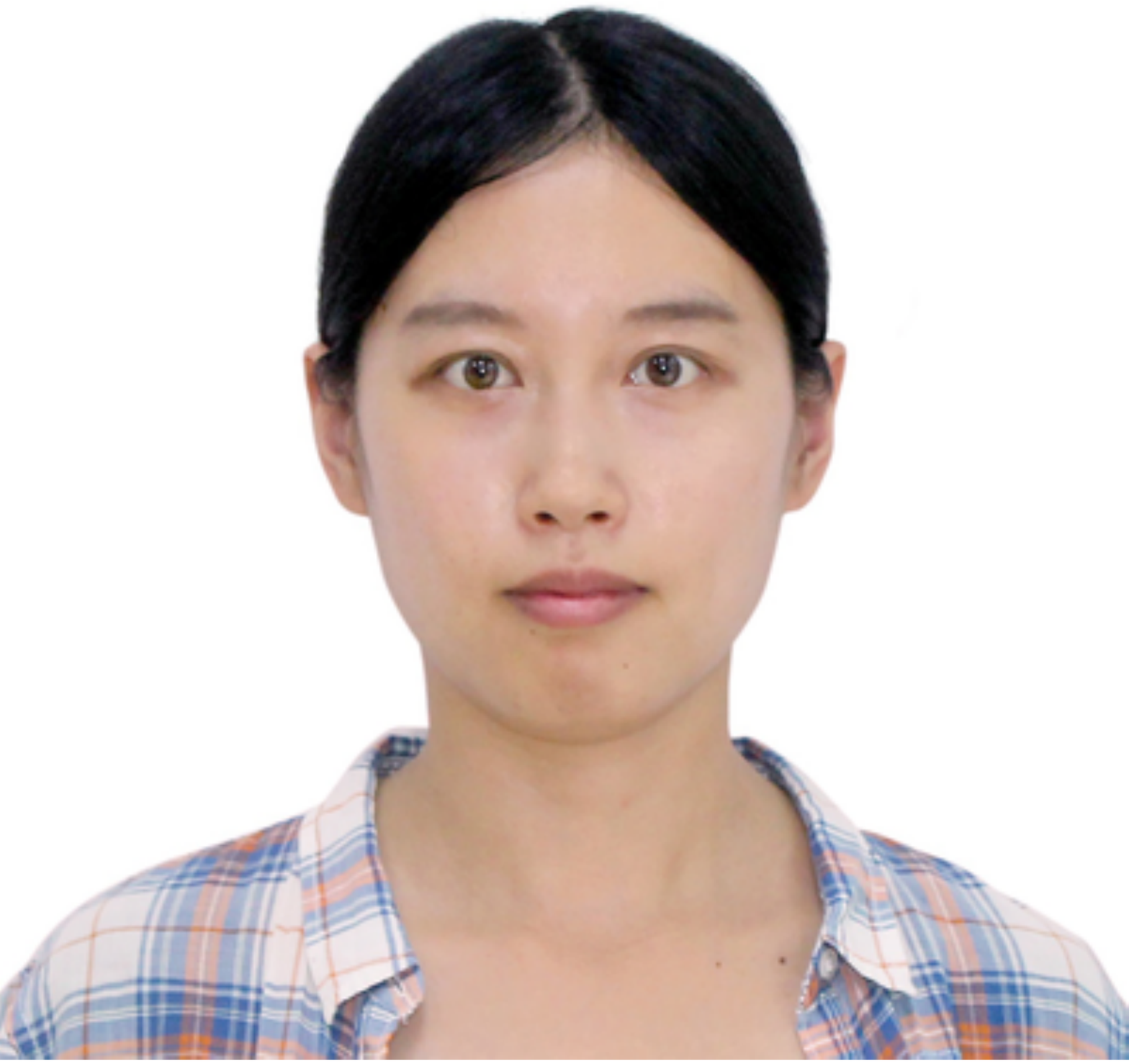}}]{Chunfeng Cui} received the Ph.D. degree in computational mathematics from Chinese Academy of Sciences, Beijing, China, in 2016 with a specialization in numerical optimization. From 2016 to 2017, she was a Postdoctoral Fellow at City University of Hong Kong, Hong Kong. In 2017 She joined the Department of Electrical and Computer Engineering at University of California Santa Barbara as a Postdoctoral Scholar.

Dr. Cui's research activities are mainly focused on the areas of tensor computing, uncertainty quantification, machine learning, and their interface.  She is the recipient of the 2019 Rising Stars in Computational and Data Sciences, 2019 Rising Stars in EECS,   the 2018 Best Paper Award of IEEE Electrical Performance of Electronic Packaging and Systems (EPEPS), and the Best Journal Paper Award of Scientia Sinica Mathematica.

\end{IEEEbiography}

\begin{IEEEbiography}
  [{\includegraphics[width=1in,height=1.25in,clip,keepaspectratio]{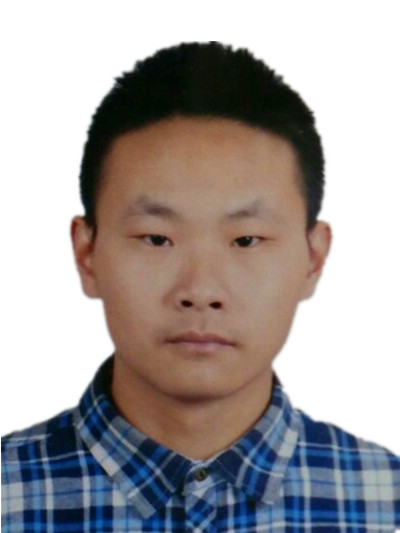}}]{Kaikai Liu} received the B.S. degree in physics in 2018 from Huazhong University of Science and Technology, Wuhan, China. In 2018 he joined the Department of Electrical and Computer Engineering at University of California Santa Barbara as a Ph.D. student.

Kaikai's research focuses on developing the novel data-driving algorithms for photonic integrated circuits design. He has been working on the tensorized Bayesian optimization and the chance-constraint optimization method.

\end{IEEEbiography}

\begin{IEEEbiography}
 [{\includegraphics[width=1in,height=1.25in,clip,keepaspectratio]{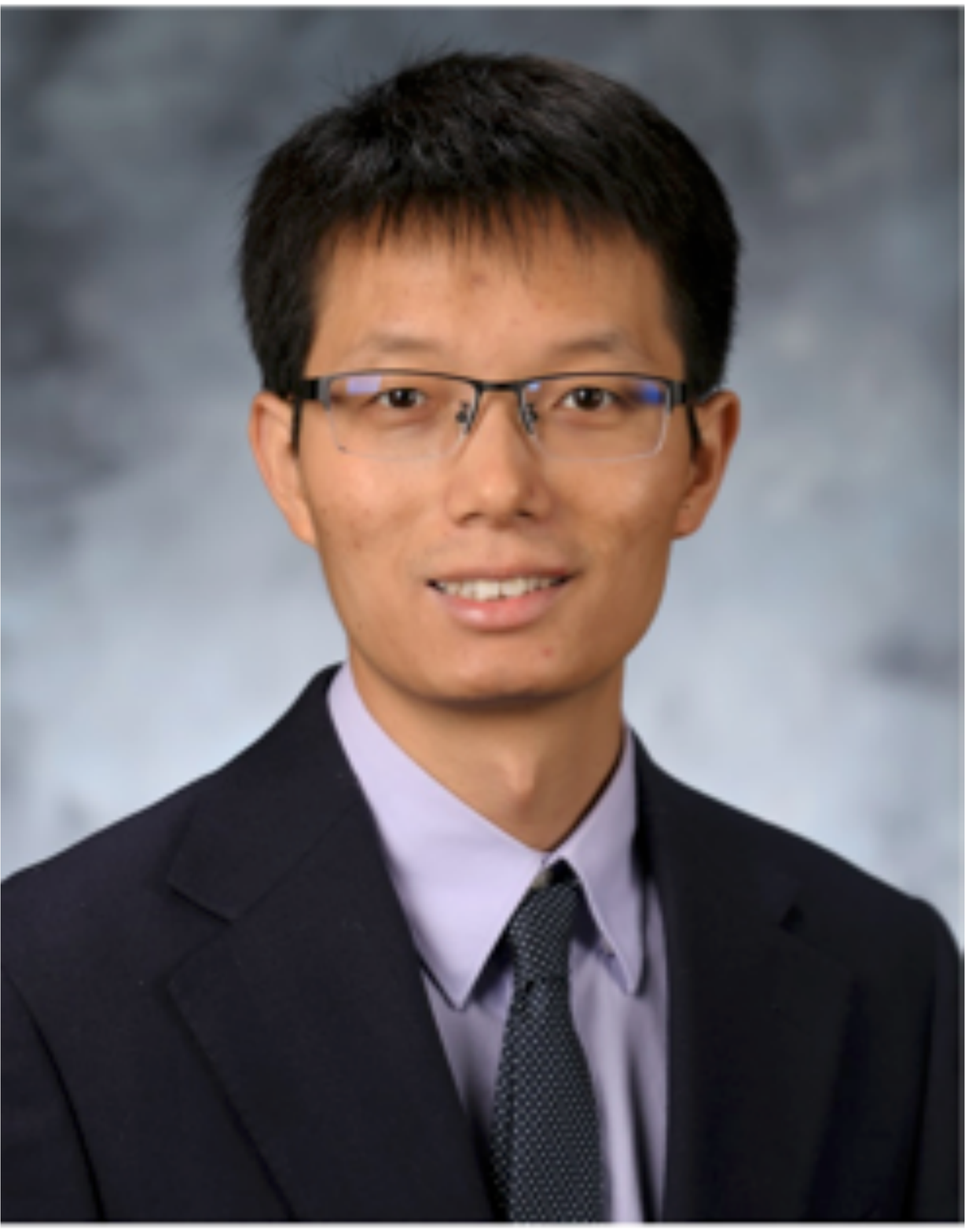}}]{Zheng Zhang} (M'15) received his Ph.D degree in Electrical Engineering and Computer Science from the Massachusetts Institute of Technology (MIT), Cambridge, MA, in 2015. He is an Assistant Professor of Electrical and Computer Engineering with the University of California at Santa Barbara (UCSB), CA. His research interests include uncertainty quantification for the design automation of multi-domain systems, and tensor methods for high-dimensional data analytics. 

Dr. Zhang received the Best Paper Award of IEEE Transactions on Computer-Aided Design of Integrated Circuits and Systems in 2014, the Best Paper Award of IEEE Transactions on Components, Packaging and Manufacturing Technology in 2018, and two Best Conference Paper Awards (IEEE EPEPS 2018 and IEEE SPI 2016). His Ph.D. dissertation was recognized by the ACM SIGDA Outstanding Ph.D. Dissertation Award in Electronic Design Automation in 2016, and by the Doctoral Dissertation Seminar Award (i.e., Best Thesis Award) from the Microsystems Technology Laboratory of MIT in 2015. He received the NSF CAREER Award in 2019. 
\end{IEEEbiography}

\end{document}